\documentclass[10pt]{amsart}
\usepackage{amsfonts,amsmath}
\usepackage{mathtools}
\usepackage{amssymb}
\usepackage{latexsym}
\usepackage{epsfig}
\usepackage{graphicx,color,graphics}
\usepackage{appendix}

\newtheorem{lemma}{Lemma}[section]

\newtheorem{theorem}[lemma]{Theorem}

\newtheorem{remark}[lemma]{Remark}

\newtheorem{teo}{Theorem}[section]

\def\sech{\operatorname{{sech}}}
\newcommand{\RR}{\mathbb{R}}

\newcommand{\NN}{\mathbb{N}}
\newcommand{\CC}{\mathbb{C}}

\def\bfD{\boldsymbol{D}}
\def\bfF{\boldsymbol{F}}
\def\bfI{\boldsymbol{I}}
\def\bfJ{\boldsymbol{J}}

\def\bfp{\boldsymbol{p}}

\def\bfu{\boldsymbol{u}}

\def\bfD{\boldsymbol{D}}
\def\bfF{\boldsymbol{F}}
\def\bfI{\boldsymbol{I}}
\def\bfId{\boldsymbol{Id}}
\def\bfJ{\boldsymbol{J}}

\def\bfp{\boldsymbol{p}}

\def\bfu{\boldsymbol{u}}

\def\O{{\mathcal O}}

\begin{document}

\title{Higher order nonlinear
Schr\"{o}dinger equation in domains with moving boundaries}

\author{Octavio V. Villagran}
\thanks{Departamento de Matemática, Universidad
de Tarapac\'a, Arica,
Chile, {\tt opverav@academicos.uta.cl; octaviovera49@gmail.com}.   The author is partially financed by projects UTA MAYOR 2022-2023, 4764-22 and UTA MAYOR 2023-2024, 4772-23.}

\author{M. Sep\'ulveda}
\thanks{DIM and CI$^2$MA, Universidad
de Concepci\'on,  Concepci\'on,
Chile, {\tt maursepu@udec.cl; masepulveda.cortes@gmail.com}.   The author is partially financed by project Fondecyt 1220869. ANID-Chile, ANID project ECOS200018, and by Centro de Modelamiento Matemático (CMM), ACE210010 and FB210005, BASAL funds for centers of excellence from ANID-Chile}

\author{R. V\'ejar Asem}
\thanks{Departamento de Matemáticas, Universidad
de la Serena, La Serena
Chile, {\tt rodrigo.vejar@userena.cl}. }

\author{Ra\'{u}l Nina Mollisaca}
\thanks{Programa Doctorado en Ciencias con Menci\'{o}n en Matem\'{a}tica, Universidad
de Tarapac\'{a}, Av. 18 de Septiembre 2222, Arica,
Chile, raulnmollisaca@gmail.com.}

\date{}

\begin{abstract}
\noindent
The initial-boundary value problem in a bounded domain with moving boundaries and nonhomogeneous boundary conditions for a higher order nonlinear Schr\"{o}dinger  (HNLS) equation is considered.
Existence and uniqueness of global weak solutions are proved as well as the stability of the solution.
Additionally, a conservative numerical method of finite differences is introduced that also verifies stability properties with respect to the $L^2$-norm, and along with proving its convergence, some interesting numerical examples are shown that illustrate the behavior of the solution.
\end{abstract}
\maketitle

{\bf Keywords}: Higher order nonlinear Schr\"{o}dinger equation, Moving boundary,
Conservative Finite Difference Method.  \\
{\it 2000 Mathematics Subject Classification}: Primary 35Q53, 35Q55;
Secondary 47J353, 35B35

\renewcommand{\theequation}{\thesection.\arabic{equation}}
\setcounter{equation}{0}

\section{Introduction}

In this work we study the existence of the weak solution as well as the asymptotic behaviour,
and numerical approximation of three-order nonlinear Schr\"{o}dinger equation (HNLS) in a bounded domain with moving boundary. Indeed, for real $\tau\geq 0,$ let $I_{\tau}$ be time-moving interval:
$
I_{\tau} = \{\xi\in\mathbb{R}:\,
\alpha(\tau) < \xi < \beta(\tau)\}$ and $Q_{\tau}$ denote bounded domain with moving boundaries:
$Q_{T} = \{(\xi,\,\tau)\in\mathbb{R}^{2}:\xi\in
I_{\tau},\,0 < \tau < T\}.$ In $Q_{T}$ we consider the higher order nonlinear Schr\"odinger equation
\begin{eqnarray}
	\label{101}
	\begin{cases}
		iv_{\tau} + \gamma v_{\xi\xi} + i\chi v_{\xi\xi\xi} = |v|^{2}v
		\quad \mbox{in } Q_{T},   \\
		v(\xi,\,0)  = v_{0}(\xi)\quad\mbox{in}\quad
		I_{0}, \\
		v(\alpha(\tau),\,\tau) =  v(\beta(\tau),\,\tau) = v_{\xi}(\alpha(\tau),\,\tau) =v_{\xi}(\beta(\tau),\,\tau)=  0, \\
		v_{\xi\xi}(\alpha(\tau),\,\tau) =  v_{\xi\xi}(\beta(\tau),\,\tau) = v_{\xi\xi\xi}(\alpha(\tau),\,\tau) =v_{\xi\xi\xi}(\beta(\tau),\,\tau)=  0, \\
	\end{cases}
\end{eqnarray}
where $\gamma$ and $\chi$ are real positive constants and $i^2 = -1.$
Moreover, we assume that there exists $\alpha_{0}$ and $\beta_{0}$ such that
\begin{eqnarray}
	\label{102}0<\alpha_{0}\leq \beta(\tau) - \alpha(\tau)\leq
	\beta_{0}<+\infty,\quad \forall\;\tau\geq 0.
\end{eqnarray}
The lateral boundary $\Sigma_{\tau}$ of $Q_{\tau}$ is given by
\begin{eqnarray}
	\label{105}\Sigma_{\tau} =
	\bigcup_{0<\tau<T}\left(\alpha(\tau)\times \{\tau\}\right)\cup
	\left(\beta(\tau)\times \{\tau\}\right).
\end{eqnarray}
The moving boundary and initial value problem
\eqref{101}
is motivated by the generalized HNLS Cauchy problem
\begin{align}
	\label{1001}
	\begin{cases}
		iv_{t} + \gamma v_{xx} + i\chi v_{xxx} + \vartheta |v|^{2}v + i\delta |v|^{2}v_{x} + i\varepsilon v^{2}\overline{v}_{x} = 0,\quad x,\ t\in \mathbb{R} \\
		v(x,\,0) = v_{0}(x)
	\end{cases}
\end{align}
where $\gamma,$ $\chi,$ $\vartheta,$ and $\delta$ are real numbers with $\chi \neq 0.$ This equation was first proposed Hasegawa {\it et al.} \cite{hasegawa} as a model for the propagation of a signal in an optical fiber (see also \cite{kodama}). The equation given above \eqref{1001} can be reduced to other well-known equations. For instance, setting
$\chi= \delta=\varepsilon = 0$ we have the classical cubic semilinear Schr\"odinger equation, that is,
\begin{align}
	\label{10001}
	\begin{cases}
		iv_{t} + \gamma v_{xx} + \vartheta |v|^{2}v  = 0,\quad x,\ t\in \mathbb{R}  \\
		v(x,\,0) = v_{0}(x).
	\end{cases}
\end{align}
If $\gamma = 1,$ $\chi=\vartheta = 0$ we have the derivative Schr\"odinger equation, that is,
\begin{align}
	\label{100001}
	\begin{cases}
		iv_{t} + v_{xx} + i\delta |v|^{2}v_{x} + i\varepsilon v^{2}\overline{v}_{x} = 0,\quad x,\ t\in \mathbb{R}  \\
		v(x,\,0) = v_{0}(x).
	\end{cases}
\end{align}
This equation \eqref{100001} is encountered in plasma physics (propagation of Alfven and radio waves) \cite{Chen}.
Letting $\gamma=\vartheta = \varepsilon =0$ the equation that arises is the
Korteweg-de Vries equation, this is,
\begin{align}
	\label{1000001}
	\begin{cases}
		v_{t} + \chi v_{xx} + \delta |v|^{2}v_{x}  = 0\quad x,\ t\in \mathbb{R}  \\
		v(x,\,0) = v_{0}(x).
	\end{cases}
\end{align}
The initial value problem for the equations \eqref{10001}, \eqref{100001} and \eqref{1000001} has been extensively studied in the last few years. See, for instance, \cite{bi01, bi02, bona} and references therein. In 1992, Laurey \cite{laurey} considered the equation \eqref{1001} and proved local well-posedness of the initial value problem associated for data in $H^{s}(\mathbb{R}),$ $s > 3/4,$ and global well-posedness in $H^{s}(\mathbb{R}),$ $s \geq 3/4.$ In 1997, Staffilani \cite{Staffilani} for  \eqref{1001} established local well-posedness for data in $H^{s}(\mathbb{R}),$ $s \geq 1/4$ improving Laurey’s result. A similar result was given in \cite{xavier} with $w(t),\,\gamma(t)$ real functions. \\ \\
Applying the Gauge transform to the HNLS equation on the real line \eqref{1001}, the $u_{xx}$ term can be removed, and the equation becomes a modified complex Korteweg-de Vries \cite{alves}. In the case of our equation with moving boundary \eqref{101}, the Gauge transform as a bi-continuous map from $L^{p}([0,\,T]: H^{s}(Q_\tau)$ to $L^{p}([0,\,T]: H^{s}(\widetilde{Q}_t)$, allows us to obtain a semilinear Airy equation with a complex nonlinear reaction term and boundary conditions in another moving boundary domain, which is not necessarily simpler than the original equation \eqref{101} (see appendix A).
Thus, this one of the motivations to continue this study for the equations \eqref{101}.
\\
\\
In this paper, we study initial-boundary value problem for the higher order nonlinear problem equation \eqref{101} in a bounded domain with moving boundaries and nonhomogeneous boundary conditions. Following the idea of
transforming the moving boundary given by Doronin and Larkin \cite{do}
(see also \cite{bi,octavio}),
we study the HNLS equation in domains with moving boundaries where now the function $v$ is a complex valued function. In our case, the aim of this work is to establish the well-posedness of the system with moving boundaries as well as the exponential decay of small solutions in asymptotically cylindrical domains accompanied by numerical results. The most common way to treat this kind of problems is to transform the domain with moving boundaries into a cylindrical one. After the transformation of the domain, we obtain a HNLS equation  type with coefficients depending on the space and time variables. Then, we study the existence, uniqueness and stability for this new problem and we give some numerical examples. For this, we consider a
Crank-Nicolson
finite difference scheme based on
conservation laws
similar to the one described in
\cite{rodrigo1, rodrigo2, rodrigo3}.
This method is originally inspired by L2 norm-preserving, stable, and convergent schemes for dispersive equations
\cite{colin, mauri} for the KdV and the KdV–Kawahara equations. Our method is different from the one presented in \cite{gear1, gear2} where the authors used Fourier transform for the x-variable. The advantage of our method is that it is better adapted to our nonperiodic boundary conditions, and to the treatment of the system with variable coefficients in the cylindrical domain given by the problem with moving boundaries. \\
\\
Our paper is organized as follows: In Section 2 we transform the moving boundary domain into the cylindrical one. In Section 3 and 4 we obtain the existence and uniqueness of solutions to the model. in Section 5 we analyze the solution of \eqref{101} using suitable multiplier techniques. Finally in Section 6 a numerical scheme is given. \\
\\
Throughout this paper $C$ is a generic constant, not necessarily the same at each occasion (it will change from line to line), which depends in an increasing way on the indicated quantities.

\renewcommand{\theequation}{\thesection.\arabic{equation}}
\setcounter{equation}{0}

\section{Preliminaries}
In this section, some notations will be introduced. Let $\Omega\subset\mathbb{R}$ open and bounded, the Lebesgue space
$L^p(\Omega)$, $1\le p\le+\infty$, with norm denoted by $\|\cdot\|_{L^p}$ will be used.
Let $m$ be a positive integer, i.e., $m\in\mathbb{N}$, the usual Sobolev spaces
$W^{m,p}(\Omega)=\{u\in L^p(\Omega)\,:\, \|D^\alpha
u\|_{L^p}<+\infty,\ \forall |\alpha|\le m\}$, with norm denoted by
$\|\cdot\|_{W^{m,p}}$ is considered. When $p=2$,  it is stablished
$H^m(\Omega):=W^{m,2}(\Omega)$  denoting the respective norm by
$\|\cdot\|_{H^m}$. Denote by $C_0(\Omega)$ the space of continuous functions with compact support in $\Omega$ and $C_0^\infty(\Omega)=C^\infty(\Omega)\cap C_0(\Omega)$ also by $H_0^1(\Omega)$ the closed of $C_0^\infty(\Omega)$ in $H^1(\Omega)$. \\
For any interval $I$ of $\mathbb{R}$ and any Banach space $X$ with the norm $\|\cdot \|_{X}$, $L^p(I:\,X)$  indicates the space of valued
functions in $X$ defined on the interval $I$ that are integrable
in the Bochner sense, and its norm will be denoted by
$\|\cdot\|_{L^p(I:X)}$. Denote $C(0,T:X)$ the space of continuous functions form $[0,T]$ into $X$ equipped with the uniform convergence norm.
\subsection{Moving Boundary Transformation}
In this section we will provide the details of the transform of the system of equations with moving boundary $Q_{T}$ into the cylindrical one. For this, we consider the application
\begin{eqnarray}
	\label{201}\left.
	\begin{array}{l}
		g:Q_{T}\rightarrow Q \\
		(\xi,\,\tau)\rightarrow \left(\frac{\xi -
			\alpha(\tau)}{\beta(\tau) - \alpha(\tau)},\,\tau\right) = (x,\,t).
	\end{array}
	\right.
\end{eqnarray}
where $Q = \{x\in(0,\,1),\ t\in (0,\,T)\}.$
Note that $g\in C^{2}$ and its inverse $g^{-1}$ is given by
\begin{eqnarray}
	\label{202}\left.
	\begin{array}{l}
		g^{-1}:Q\rightarrow Q_{T} \\
		(x,\,t)\rightarrow (\xi,\,\tau)=\left(\alpha(t) + \left(\beta(t) -
		\alpha(t)\right)x,\,t\right),
	\end{array}
	\right.
\end{eqnarray}
and it is also $C^{2}.$ Let ${\bf p}(\tau) = \frac{1}{\beta(\tau) -
	\alpha(\tau)}.$ Denoting
\begin{eqnarray}
	\label{203}u(x,\,t) = v\circ g^{-1}(x,\,t) = v(\xi,\,\tau)
\end{eqnarray}
and taking into account the change of variable \eqref{201}, we have
\begin{align*}
	v_{\tau} &= \frac{\partial
		u}{\partial x}\frac{\partial x}{\partial\tau} + \frac{\partial
		u}{\partial t}\frac{\partial t}{\partial\tau} = \frac{\partial
		u}{\partial x}\frac{\partial x}{\partial\tau} + \frac{\partial
		u}{\partial t} \\
	&= u_{x}\left[-\frac{\alpha'(\tau)}{\beta(\tau)-\alpha(\tau)}-\frac{(\xi-\alpha(\tau))(\beta'(\tau)-\alpha'(\tau))}{(\beta(\tau)-\alpha(\tau))^2}\right]+
	u_{t}\\
	&= u_{x}\left[\alpha'(\tau){\bf p}(\tau)+\frac{(\xi-\alpha(\tau))}{(\beta(\tau)-\alpha(\tau))}(\ln {\bf p}(\tau))'\right]+
	u_{t}\\
	&= u_{x}\left[\alpha'(\tau){\bf p}(\tau)+x(\ln {\bf p}(\tau))'\right]+
	u_{t}
\end{align*}
where the notation $\ ^{\pmb{\color{black}{{'}}}} = \frac{d}{d\tau}.$
Let $L(\xi,\,\tau)=\alpha'(\tau){\bf p}(\tau)+x(\ln {\bf p}(\tau))'$, then
\begin{align}\label{204}
	v_{\tau}(\xi,\, \tau) &= u_{x}(g(\xi,\ \tau))L(\xi,\ \tau)+
	u_{t}(g(\xi,\, \tau))
\end{align}
for all $(\xi,\,\tau)\in Q_T$. Similary
\begin{align*}
	&v_{\xi}(\xi,\,\tau)= u_x(g(\xi,\,\tau)){\bf p}(\tau),\hspace{1cm}
	v_{\xi\xi}(\xi,\,\tau) = u_{xx}(g(\xi,\,\tau))[{\bf p}(t)]^{2}\\
	&v_{\xi\xi\xi}(\xi,\,\tau) = u_{xxx}(g(\xi,\,\tau))[{\bf p}(\tau)]^{3}
\end{align*}
for all $(\xi,\,\tau)\in Q_T$.
Hence, into \eqref{101}$_{1}$ we have
\begin{align}
	&u_{t}(g(\xi,\,\tau)) -i\gamma u_{xx}(g(\xi,\,\tau))[{\bf p}(\tau)]^2 +\chi [{\bf p}(t)]^{3}u_{xxx}(g(\xi,\,\tau))+ L(\xi,\,\tau)u_{x}(g(\xi,\,\tau))\nonumber\\
	& = -i|u(g(\xi,\,\tau))|^{2}u(g(\xi,\,\tau)).
\end{align}
How $(x,\,t)=g(\xi,\,\tau)$ and $L(\xi,\,\tau)=-\alpha'(\tau){\bf p}(\tau)+x[\ln{\bf p}(\tau)]'=-\alpha'(t){\bf p}(t)+x[\ln{\bf p}(t)]'$, we denote $L(x,\,t)=-\alpha'(t){\bf p}(t)+x[\ln{\bf p}(t)]'$, then
\begin{align}
	\label{207}u_{t} - i\gamma u_{xx}[{\bf p}(t)]^2 +\chi [{\bf p}(t)]^{3}u_{xxx}+ L(x,\,t) u_{x} +i|u|^{2}u=0
\end{align}
for all $(x,\,t)\in Q.$\\
We denote by $I(x,\,t)=i L(x,\,t)$, $\mathfrak{B}(t)=[{\bf p}(t)]^{3}$ and $\mathfrak{C}(t)=[{\bf p}(t)]^2$, then
\begin{align}
	\label{eq57}i u_{t} +i\chi \mathfrak{B}(t)u_{xxx}+\gamma\mathfrak{C}(t) u_{xx} + I(x,\,t) u_{x} -|u|^{2}u=0
\end{align}
On the other hand, on the boundary, the variable $x$ becomes
\begin{eqnarray}
	\label{208} x = \frac{\xi - \alpha(\tau)}{\beta(\tau) -
		\alpha(\tau)} = \left\lbrace
	\begin{array}{l}
		0,\quad \mbox{if}\quad \xi = \alpha(\tau), \\
		\\
		1,\quad \mbox{if}\quad \xi = \beta(\tau),
	\end{array}
	\qquad \Longrightarrow \quad 0<x<1,
	\right.
\end{eqnarray}
that is, $\xi = \left(\beta(0) -
\alpha(0)\right)x + \alpha(0)$
then the initial condition is
\begin{eqnarray}
	\label{209}u(x,\,0) = (v\circ g^{-1})(x,0)=u_0(\alpha(0)+(\beta(0)-\alpha(0))x).
\end{eqnarray}
and the boundary conditions remain homogeneous.

\subsection{Main Result}

We consider the system of PDE
\begin{align}\label{eq47}
	&iu_{t} + i\chi\mathfrak{B}(t)u_{xxx} + \gamma  \mathfrak{C}(t)u_{xx} + I(x,\,t)u_{x} - |u|^{2}u = 0
\end{align}
$(x,\,t)\in Q=(0,\,1)\times (0,\,T)$ and with initial and boundary conditions
\begin{align}
	\label{eq48}& u(x,\,0) = u_{0}(x)\\
	\label{eq49}& u(0,\,t) = u_x(0,\,t) = u_{xx}(0,\,t) = u_{xxx}(0,\,t) = 0\\
	\label{eq50}& u(1,\,t) = u_x(1,\,t) = u_{xx}(1,\,t) = u_{xxx}(1,\,t) = 0.
\end{align}
Our main result is as follows:
\begin{enumerate}
	\item Let $u_0(x)\in H^4(0,1)$, there is one classic solution of \eqref{eq47}-\eqref{eq50} (or \eqref{101}).
	\item There is only one solution of \eqref{eq47}-\eqref{eq50} (or \eqref{101}).
	\item If $u$ is solution of \eqref{eq47}-\eqref{eq50}, then
	\begin{align}
		\|u\|_{L^2(0,1)}^2 \leq Ce^{-t}
	\end{align}
	where $C$ depends only $u_0$.
\end{enumerate}

\renewcommand{\theequation}{\thesection.\arabic{equation}}
\setcounter{equation}{0}

\section{Existence of solutions}
In the previous Section we have transformed a domain with moving boundaries into cylindrical one. Now we will apply the Faedo-Galerkin method in the transformed mixed problem.\\
For this we study the regularized problem with a spatial basis and we prove the existence of global solutions for this regularization. This regularization procedure is classical adding a dissipative term of fourth order for a KdV type equation (see \cite{lions}). In our case, we consider a basis and a regularization very close to the one used in \cite{do, la}. Let $\varepsilon > 0.$ For the unknown functions $u^{\varepsilon}(x,\,t)$ we consider the following regularized problem (hereafter we omit the index $\varepsilon$ in calculations):
\begin{align}
	\label{301}&iu_{t} + i\chi\mathfrak{B}(t)u_{xxx} + \gamma  \mathfrak{C}(t)u_{xx} + I(x,\,t)u_{x} + i\varepsilon \mathfrak{B}(t)u_{xxxx} - |u|^{2}u = 0
\end{align}
with $Q =\in (0,\,1)\times (0,\,T)$ and initial and boundary conditions
\begin{align}
	& \label{eq52}u(x,\,0) = u_{0}(x) \\
	& \label{eq53}u(0,\,t) = u_x(0,\,t) = u_{xx}(0,\,t) = u_{xxx}(0,\,t) = 0\\
	& \label{eq54}u(1,\,t) = u_x(1,\,t) = u_{xx}(1,\,t) = u_{xxx}(1,\,t) = 0.
\end{align}
We assume that $3\chi-4\varepsilon>0$.
\begin{theorem} Let $u_0\in H^4(0,1)$, there is one classic solution of \eqref{eq47}-\eqref{eq50} (or \eqref{101}).
\end{theorem}
\proof  Let $\{g_{j}(x)\}_{j\in\mathbb{N}}$ the Hilbertian base of $L^{2}(0,\,1)$ solution of the eigenvalue problem
\begin{align}
	\label{305}&g_{j}(0) = g_{j}'(0)= g_j''(0)=g_j'''(0)=0, \qquad j\in\mathbb{N}\\
	\label{eq55}& g_j(1) = g_{j}'(1)= \varepsilon g_{j}''(1) = g_j'''(1)= 0,\qquad j\in\mathbb{N}.
\end{align}
We construct the approximation
\begin{align}
	\label{306}u^{N}(x,\,t) = \sum_{j = 1}^{N}
	h_{j}^{N}(t)g_{j}(x),
\end{align}
where $h_{j}^{N}(t)g_{j}(x)$ is solution for the nonlinear Schr\"{o}dinger equation
\begin{align}\label{307}
	&i\langle u_{t}^{N},\,g_{j}\rangle + i\chi\mathfrak{B}(t)\langle u_{xxx}^{N},\,g_{j}\rangle + \gamma \mathfrak{C}(t)\langle u_{xx}^{N},\,g_{j}\rangle + \langle I(x,\,t)u_{x}^{N},\,g_{j}\rangle \nonumber \\
	&\qquad + i\varepsilon\mathfrak{B}(t)\langle u_{xxxx}^{N},\,g_{j}\rangle
	- \langle |u^{N}|^{2}u^{N},\,g_{j}\rangle =  0
\end{align}
for $j = 1,\ldots,\,N$ and subject to initial data
\begin{eqnarray}
	\label{308}h_{j}^{N}(0) = \langle u_{0},\,g_{j}\rangle=0
\end{eqnarray}
such that $\langle\cdot,\,\cdot\rangle\equiv  \langle\cdot,\,\cdot\rangle_{[L^{2}(0,\,1)]^{2}}.$ Due to orthonormality of $\{g_{j}\}_{j\in\mathbb{N}},$ the coupled system \eqref{307} can be written in the norm form. Moreover, according to the standard existence theory for ordinary differential equations, a unique solution exists on some time-interval.

From the nonlinear multipliers we are able to pass to the limits in the linear term of \eqref{307}. The nonlinear terms need a careful analysis. For that, we will need to use the Lions–Aubin compactness theorem (see \cite{au, lions}), as we will see in the next subsections. \\
{\it First estimate.}
Setting $g_{j} = e^{x}\overline{u}^{N}$ in \eqref{307} we have
\begin{align*}
	&i\langle u_{t}^{N}\,e^{x},\,\overline{u}^{N}\rangle + i\chi \mathfrak{B}(t)\langle u_{xxx}^{N}e^{x},\,\overline{u}^{N}\rangle + i\varepsilon \mathfrak{B}(t)\langle u_{xxxx}^{N}e^{x},\,\overline{u}^{N}\rangle \nonumber \\
	&\qquad + \langle I(x,\,t)u_{x}^{N}e^{x},\,\overline{u}^{N}\rangle + \gamma  \mathfrak{C}(t)\langle u_{xx}^{N}e^{x},\,\overline{u}^{N}\rangle - \langle |u^{N}(t)|^{2}u^{N}e^{x},\,\overline{u}^{N}\rangle =  0.
\end{align*}
Note that integrating by parts and using \eqref{eq52}-\eqref{eq54} we have
\begin{align*}
	\langle u_{xxx}^{N}e^{x},\,\overline{u}^{N}\rangle
	&=\int_0^1 u_{x}^{N}e^{x}\overline{u}^{N}\ dx+\int_0^1 u_{x}^{N}e^{x}\overline{u}_x^{N}\ dx-\int_0^1 u_{xx}^{N}e^{x}\overline{u}_x^{N}\ dx.\\
	\langle u_{xxxx}^{N}e^{x},\,\overline{u}^{N}\rangle
	&=-\int_0^1 u_{x}^{N}e^{x}\overline{u}^{N}\ dx-\int_0^1 u_{x}^{N}e^{x}\overline{u}_x^{N}\ dx+2\int_0^1 u_{xx}^{N}e^{x}\overline{u}_x^{N}\ dx\\
	&\ \ +\int_0^1 e^{x}|u_{xx}^{N}|^2\ dx.\\
	\langle u_{xx}^{N}e^{x},\,\overline{u}^{N}\rangle&
 =\int_0^1 u_{xx}^{N}e^{x}\overline{u}^{N}\ dx
	=-\int_0^1 u_{x}^{N}e^{x}\overline{u}^{N}\ dx-\int_0^1 e^{x}|u_x^{N}|^2\ dx
\end{align*}
Therefore, replacing these three expressions in the previous equation we obtain
\begin{align*}
	&i\int_{0}^{1}e^{x}u_{t}^{N}\overline{u}^{N}dx + i\chi\mathfrak{B}(t)\int_{0}^{1}e^{x}\overline{u}^{N}u_{x}^{N}dx + i\chi\mathfrak{B}(t)\int_{0}^{1}e^{x}|u_{x}^{N}|^{2}dx \nonumber \\
	& -i\chi\mathfrak{B}(t)\int_{0}^{1}e^{x}\overline{u}_{x}^{N}u_{xx}^{N}dx
	- \gamma  \mathfrak{C}(t)\int_{0}^{1}e^{x}\overline{u}^{N}u_{x}^{N}dx
	- \gamma  \mathfrak{C}(t)\int_{0}^{1}e^{x}|u_{x}^{N}|^{2}dx\nonumber\\
	&+\int_{0}^{1}I(x,\,t)e^{x}\overline{u}^{N}u_{x}^{N}dx - i\varepsilon\mathfrak{B}(t)\int_{0}^{1}e^{x}\overline{u}^{N}u_{x}^{N}dx-i\varepsilon\mathfrak{B}(t)\int_{0}^{1}e^{x}|u_{x}^{N}|^{2}dx\nonumber\\
	&+2i\varepsilon\mathfrak{B}(t)\int_{0}^{1}e^{x}\overline{u}_{x}^{N}u_{xx}^{N}dx
	+i\varepsilon\mathfrak{B}(t)\int_{0}^{1}e^{x}|u_{xx}^{N}|^{2}dx- \int_{0}^{1}e^{x}|u^{N}|^{4}dx\\
	& = 0.
\end{align*}
Considering then the imaginary part of this expression and applying integration by parts, we deduce
\begin{align}\label{eq16}
	&\frac{d}{dt}\langle e^{x},\,|u^{N}|^{2}\rangle - (\chi - \varepsilon)\mathfrak{B}(t) \langle e^{x},\,|u^{N}|^{2}\rangle
	+\mathfrak{B}(t)(3\chi-4\varepsilon) \langle e^{x},\,|u_x^{N}|^{2}\rangle+ 2\varepsilon\mathfrak{B}(t)\langle e^{x},\,|u_{xx}^{N}|^{2}\rangle\nonumber\\
	&=\ 2\gamma \mathfrak{C}(t) Im\langle e^{x}\overline{u}^{N},\,u_{x}^{N}\rangle
	- 2Im\langle I(x,\,t)e^{x}\overline{u}^{N},\,u_{x}^{N}\rangle.
\end{align}
We are going to bounded the right side of \eqref{eq16}
\begin{align*}
	&2\gamma \mathfrak{C}(t) Im\langle e^{x}\overline{u}^{N},\,u_{x}^{N}\rangle
	\leq C\|e^x\overline{u}^N\|_{L^2(0,\,1)}\|\overline{u}_x^Ne^x\|_{L^2(0,\,1)}
	\leq \varepsilon_1\|\overline{u}_x^Ne^x\|_{L^2(0,\,1)}^2+
 +C(\varepsilon_1)\|e^x\overline{u}^N\|_{L^2(0,\,1)}^2\\
	&\leq \varepsilon_1e\int_0^1e^x|u^N|^2\ dx+C(\varepsilon_1)e\int_0^1e^x|u_x^N|^2\ dx.\\
	&- 2Im\langle I(x,\,t)e^{x}\overline{u}^{N},\,u_{x}^{N}\rangle
	\leq 2\|I(x,\,t)e^xu_x^N\|_{L^2(0,\,1)}\|\overline{u}^Ne^x\|_{L^2(0,\,1)}
	\leq\varepsilon_2\|I(x,\,t)e^xu_x^N\|_{L^2(0,\,1)}^2\\
	&+C(\varepsilon_2)\|\overline{u}^Ne^x\|_{L^2(0,\,1)}^2
	\leq \varepsilon_2\sup_{x\in(0,1)}|I(x,\,t)|^2e\int_0^1e^x|u_x^N|^2\ dx+C(\varepsilon_2)e\int_0^1e^x|u^N|^2\ dx\\
	&\leq \varepsilon_3\int_0^1e^x|u_x^N|^2\ dx+C(\varepsilon_3)\int_0^1e^x|u^N|^2\ dx.
\end{align*}
Then
\begin{align}\label{eq12}
	&\frac{d}{dt}\langle e^{x},\,|u^{N}|^{2}\rangle +(3\chi-4\varepsilon-\varepsilon_1-\varepsilon_3)\ \langle e^{x},\,|u_x^{N}|^{2}\rangle+2\varepsilon\mathfrak{B}(t)\langle e^{x},\,|u_{xx}^{N}|^{2}\rangle\nonumber\\
	&\leq (|\chi-\varepsilon|+C(\varepsilon_1)e+C(\varepsilon_3))\langle e^{x},\,|u^{N}|^{2}\rangle.
\end{align}
For $\varepsilon_1$ and $\varepsilon_3$ suitable and integrating over $(0,t)$ with $t\leq T$ it follows that
\begin{align*}
	&\frac{d}{dt}\langle e^{x},\,|u^{N}|^{2}\rangle\leq\langle e^{x},\,|u^{N}(\cdot,\,0)|^{2}\rangle
	+C\int_0^t\langle e^{x},\,|u^{N}|^{2}\rangle\ ds\\
	&\leq e\|u_0^N\|_{L^2(0,\,1)}^2+C\int_0^t\langle e^{x},\,|u^{N}|^{2}\rangle\ ds.
\end{align*}
By Gronwall inequality we have
\begin{align}\label{eq22}
	\|u^N\|_{L^2(0,\,1)}^2\leq \langle e^{x},\,|u^{N}|^{2}\rangle\leq C\exp(Ct)
\end{align}
for all $t\in (0,T)$. Therefore
\begin{align}\label{eq23}
	\textrm{ $(u^N)$ is bounded in $L^\infty(0,T:L^2(0,1))\hookrightarrow L^2(0,T:L^2(0,1)) $}
\end{align}
Note that of \eqref{eq12} we have
\begin{align}\label{eq13}
	&\langle e^{x},\,|u^{N}|^{2}\rangle +(3\chi-4\varepsilon-\varepsilon_1-\varepsilon_3)\int_0^t \langle e^{x},\,|u_x^{N}|^{2}\rangle\ ds
	+2\varepsilon\int_0^t\mathfrak{B}(s)\langle e^{x},\,|u_{xx}^{N}|^{2}\rangle\ ds\nonumber\\
	&\leq e\|u_0\|_{L^2(0,\,1)}^2+C
	\leq C\|u_0\|_{L^2(0,\,1)}^2.
\end{align}
Therefore
\begin{align}\label{eq24}
	\textrm{ $(u_x^N)$   is bounded in $L^2(0,T:L^2(0,1)) $}
\end{align}
\begin{remark}
	\begin{enumerate}
		\item If $\varepsilon=0$, for $\varepsilon_1$ and $\varepsilon_3$ suitable, we have
		\eqref{eq22}-\eqref{eq23} and $\eqref{eq13}-\eqref{eq24}$, this is the constant $C$ independent of $\varepsilon$.
		\item Note that of \eqref{eq13} we have
		$$\textrm{ $(u_{xx}^N)$   is bounded in $L^2(0,T:L^2(0,1)) $}$$
		but the estimative depend of $\varepsilon$. \\
		This is a difficulty, we will improve the estimative for $u_{xx}^N$.
	\end{enumerate}
\end{remark}

\medskip

\noindent
{\it Second estimate.}
Setting $g_{j} = \overline{u}_{xxxx}^{N}$ in \eqref{307} we have
\begin{align}
	&i\langle u_{t}^{N},\,\overline{u}_{xxxx}^{N}\rangle + i\chi \mathfrak{B}(t)\langle u_{xxx}^{N},\,\overline{u}_{xxxx}^{N}\rangle + \gamma  \mathfrak{C}(t)\langle u_{xx}^{N},\,\overline{u}_{xxxx}^{N}\rangle \nonumber \\
	\label{404}&+ \langle I(x,\,t)u_{x}^{N},\,\overline{u}_{xxxx}^{N}\rangle
	+ i\varepsilon \mathfrak{B}(t)\langle u_{xxxx}^{N},\,\overline{u}_{xxxx}^{N}\rangle - \langle |u^{N}|^{2}u^{N},\,\overline{u}_{xxxx}^{N}\rangle =  0.
\end{align}
Then, multiplying by $\chi$ we have
\begin{align*}
	&i\chi\langle u_{t}^{N},\,\overline{u}_{xxxx}^{N}\rangle + i\chi^{2}\mathfrak{B}(t)\langle u_{xxx}^{N},\,\overline{u}_{xxxx}^{N}\rangle + \gamma\chi \mathfrak{C}(t)\langle u_{xx}^{N},\,\overline{u}_{xxxx}^{N}\rangle \nonumber \\
	&+ \chi\langle I(x,\,t)u_{x}^{N},\,\overline{u}_{xxxx}^{N}\rangle
	+ i\chi\varepsilon \mathfrak{B}(t)\|u_{xxxx}^{N}\|_{L^{2}(0,\,1)}^{2} - \chi\langle |u^{N}|^{2}u^{N},\,\overline{u}_{xxxx}^{N}\rangle =  0.
\end{align*}
Note that
\begin{align*}
	\langle u_{t}^{N},\,\overline{u}_{xxxx}^{N}\rangle =\int_0^1u_{txx}^N\overline{u}_{xx}^Ndx,
\end{align*}
then we obtain
\begin{align*}
	&i\chi\int_{0}^{1}u_{xxt}^{N}\overline{u}_{xx}^{N}dx
	+i\chi^2\mathfrak{B}(t)\int_{0}^{1}u_{xxx}^{N}\overline{u}_{xxxx}^{N}dx
	+\gamma\chi \mathfrak{C}(t)\int_0^1 u_{xx}^{N}\overline{u}_{xxxx}^{N}dx \\
	&+ \chi\int_0^1I(x,\,t)u_{x}^{N}\overline{u}_{xxxx}^{N}dx
	+i\chi\varepsilon\mathfrak{B}(t)\|u_{xxxx}^N\|_{L^2(0,1)}^2 - \chi\langle |u^{N}|^{2}u^{N},\,\overline{u}_{xxxx}^{N}\rangle\nonumber \\
	& = 0.
\end{align*}
Applying Conjugate
\begin{align*}
	&-i\chi\int_{0}^{1}\overline{u}_{xxt}^{N}u_{xx}^{N}dx
	-i\chi^2\mathfrak{B}(t)\int_{0}^{1}\overline{u}_{xxx}^{N}u_{xxxx}^{N}dx
	+\gamma\chi \mathfrak{C}(t)\int_0^1 \overline{u}_{xx}^{N}u_{xxxx}^{N}dx  \\
	&+ \chi\overline{\int_0^1I(x,\,t)u_{x}^{N}\overline{u}_{xxxx}^{N}dx}
	-i\chi\varepsilon\mathfrak{B}(t)\|u_{xxxx}^N\|_{L^2(0,1)}^2
	- \chi\overline{\langle |u^{N}|^{2}u^{N},\,\overline{u}_{xxxx}^{N}\rangle}\\
	& =0.
\end{align*}
Subtracting the above equations it follows that
\begin{align*}
	&i\chi\frac{d}{dt}\|u_{xx}^{N}\|_{L^{2}(0,\,1)}^{2}+i\chi^2\mathfrak{B}(t)\left(\int_0^1u_{xxx}^{N}\overline{u}_{xxxx}^{N}dx
	+\int_0^1\overline{u}_{xxx}^{N}u_{xxxx}^{N}dx\right)\\
	&+\gamma\chi\mathfrak{C}(t)\left(\int_0^1u_{xx}^{N}\overline{u}_{xxxx}^{N}dx-\int_0^1\overline{u}_{xx}^{N}u_{xxxx}^{N}dx\right)\\
	&+\chi \int_0^1I(x,t)u_x^N\overline{u}_{xxxx}^Ndx-\chi\overline{\int_0^1I(x,t)u_x^N\overline{u}_{xxxx}^Ndx}\\
	&+2i\chi\mathfrak{B}(t)\|u_{xxxx}^N\|_{L^2(0,1)}^2-\chi(\langle |u^{N}|^{2}u^{N},\,\overline{u}_{xxxx}^{N}\rangle-\overline{\langle |u^{N}|^{2}u^{N},\,\overline{u}_{xxxx}^{N}\rangle})= 0.
\end{align*}
Then
\begin{align*}
	&i\chi\frac{d}{dt}\|u_{xx}^{N}\|_{L^{2}(0,\,1)}^{2}
	+\gamma\chi\mathfrak{C}(t)2iIm\left(\int_0^1u_{xx}^{N}\overline{u}_{xxxx}^{N}dx\right)
	+2\chi i Im \left(\int_0^1I(x,t)u_x^N\overline{u}_{xxxx}^Ndx\right)\\
	&+2i\chi\mathfrak{B}(t)\|u_{xxxx}^N\|_{L^2(0,1)}^2-\chi2iIm\left(\langle |u^{N}|^{2}u^{N},\,\overline{u}_{xxxx}^{N}\rangle\right)=0
\end{align*}
or
\begin{align}\label{eq21}
	& \chi\frac{d}{dt}\|u_{xx}^{N}\|_{L^{2}(0,\,1)}^{2}
	+2\chi\mathfrak{B}(t)\|u_{xxxx}^{N}\|_{L^{2}(0,\,1)}^{2}\nonumber\\
	&=-2\gamma\chi \mathfrak{C}(t)Im\left(\langle u_{xx}^{N},\,\overline{u}_{xxxx}^{N}\rangle\right)
	-2\chi Im\left\langle I(x,\,t)u_x^N,\,\overline{u}_{xxxx}^N\rangle\right)\nonumber\\
	&+2\chi Im\left\langle |u^N|^2u^N,\,\overline{u}_{xxxx}^N\rangle\right).
\end{align}
We are going to bounded the right side of \eqref{eq21}
\begin{align*}
	&|\langle I(x,\,t)u_x^N,\,\overline{u}_{xxxx}^N\rangle|\leq \varepsilon_1\|u_{xxxx}^{N}\|_{L^2(0,1)}^2+C(\varepsilon_1)\|I(x,\,t)\|_{L^\infty(0,T:L^2(0,1))}^2\|u_x^N\|_{L^2(0,1)}^2\\
	&|\langle |u_N|^2u^N,\,\overline{u}_{xxxx}^N\rangle|\leq \varepsilon_2\|u_{xxxx}^{N}\|_{L^2(0,1)}^2+C(\varepsilon_2)\|u^N\|_{L^\infty(0,T:L^2(0,1))}^2\|u^N\|_{L^2(0,1)}^2\\
	&|\langle u_{xx}^N,\,\overline{u}_{xxxx}^N\rangle|\leq
 \varepsilon_3\|u_{xxxx}^N\|_{L^2(0,1)}^2+
 C(\varepsilon_3)\|u_{xx}^{N}\|_{L^2(0,1)}^2+.
\end{align*}
Integrating 0 to $t$ with $t\in(0,T)$
\begin{align}\label{406}
	&\chi\|u_{xx}^{N}\|_{L^{2}(0,\,1)}^{2}
	+ 2\chi\int_0^t\mathfrak{B}(s)\|u_{xxxx}^{N}\|_{L^{2}(0,\,1)}^{2}ds \nonumber \\
	&\leq \chi\|u_{xx}^N(\cdot,\,0)\|_{L^{2}(0,\,1)}^{2}
	+2\gamma\chi C(\varepsilon_3)\int_0^t\mathfrak{C}(t)\|u_{xx}^N\|_{L^2(0,\,1)}^2ds+2\gamma\chi \varepsilon_3\int_0^t\mathfrak{C}(t)\|u_{xxxx}^N\|_{L^2(0,\,1)}^2ds\nonumber\\
	&+2\chi\varepsilon_1\int_0^t\|u_{xxxx}^N\|_{L^2(0,\,1)}^2ds+2\chi C(\varepsilon_1)\|u^N\|_{L^\infty(0,T:L^2(0,\,1))}^2\int_0^t\|u_{x}^N\|_{L^2(0,\,1)}^2ds\nonumber\\
	&+2\chi\varepsilon_2\int_0^t\|u_{xxxx}^N\|_{L^2(0,\,1)}^2ds+2\chi C(\varepsilon_2)\|u^N\|_{L^\infty(0,T:L^2(0,\,1))}^2\int_0^t\|u^N\|_{L^2(0,\,1)}^2ds.
\end{align}
Then
\begin{align}\label{eq4}
	&\chi\|u_{xx}^{N}\|_{L^{2}(0,\,1)}^{2} +C_2\int_{0}^{t}\|u_{xxxx}^{N}\|_{L^2(0,\,1)}^{2}ds
	\leq \|u_0^N\|_{H^2(0,\,1)}^2\nonumber \\
	& \ \
	+\chi\int_0^t\|u_{xx}^{N}\|_{L^{2}(0,\,1)}^{2}ds
\end{align}
considering $\ \varepsilon_1, \varepsilon_2,\ \varepsilon_3 $ suitable we have $C_2>0$. By Gronwall inequality
\begin{align}\label{eq25}
	\chi\|u_{xx}^{N}\|_{L^{2}(0,\,1)}^{2}\leq C \exp(T)
\end{align}
for all $t\in(0,T)$, where $C$ independent of $\epsilon$,  then
\begin{align}\label{eq26}
	(u_{xx}^N)\  \textrm{is bounded in} \  L^\infty(0,T:L^2(0,1))\hookrightarrow L^2(0,T:L^2(0,1)).
\end{align}
Finally of \eqref{eq4} we obtain
\begin{align}\label{eq20}
	&\chi\|u_{xx}^{N}\|_{L^{2}(0,\,1)}^{2} +C_2\int_{0}^{T}\|u_{xxxx}^{N}\|_{L^2(0,\,1)}^{2}dt  \nonumber\\
	&\leq C\left(\|u_0\|_{H^2(0,\,1)}^2+\|J(x,\,s)\|_{L^2(0,T:L^2(0,1))}^2\right).
\end{align}
where $C$ independent of $\varepsilon$. Here we obtain that
\begin{align}\label{eq27}
	(u_{xxxx}^N)\  \textrm{is bounded in} \  L^2(0,T:L^2(0,1))
\end{align}
{\it Third estimate.}
Now from \eqref{307} we have
\begin{align}
	&i\eta(t)\langle u_{t}^{N},\,g_{j}\rangle + i\chi\langle u_{xxx}^{N},\,g_{j}\rangle + \gamma \mu(t)\langle u_{xx}^{N},\,g_{j}\rangle + \langle A(x,\,t)u_{x}^{N},\,g_{j}\rangle \nonumber \\
	\label{407}&+ i\varepsilon\langle u_{xxxx}^{N},\,g_{j}\rangle
	- \eta(t)\langle |u^{N}|^{2}u^{N},\,g_{j}\rangle =  0
\end{align}
where
\begin{align*}
	& \eta(t) := \mathfrak{B}(t)= [\beta(t) - \alpha(t)]^{3},\quad \mu(t) := [p(t)]^{-1}= \beta(t) - \alpha(t) \\
	& A(x,\,t) := [p(t)]^{-3}I(x,\,t)= [\beta(t) - \alpha(t)]^{3}I(x,\,t).
\end{align*}
Differentiating \eqref{407} with respect to $t > 0$ and we  multiply the result by $e^{x}\overline{u}_{t}^{N}$
\begin{align*}
	&i\eta'(t)\langle e^{x},\,|u_{t}^{N}|^{2}\rangle + i\eta(t)\langle e^{x},\,\overline{u}_{t}^{N}u_{tt}^{N}\rangle
	+i\chi\langle u_{xxxt}^{N},\,e^x\overline{u}_{t}^{N}\rangle
	+ \gamma\mu'(t) \langle u_{xx}^{N},\,e^x\overline{u}_{t}^{N} \rangle\\
	&+ \gamma \mu(t) \langle u_{xxt}^{N},\,e^x\overline{u}_{t}^{N}\rangle+ \langle A_{t}u_{x}^{N},\,e^{x}\overline{u}_{t}^{N}\rangle + \langle Au_{xt}^{N},\,e^x\overline{u}_{t}^{N}\rangle
	+ i\varepsilon\langle u_{xxxxt}^{N},\,e^{x}\overline{u}_{t}^{N}\rangle  \\
	&-\eta'(t)\langle |u^{N}|^{2}u^{N},\,e^x\overline{u}_{t}^{N}\rangle- \eta(t)\langle (|u^{N}|^{2})_{t}u^{N},\,e^x\overline{u}_{t}^{N}\rangle
	- \eta(t)\langle e^x,\,|u^{N}|^{4}\rangle=0.
\end{align*}
Applying Conjugate
\begin{align*}
	&-i\eta'(t)\langle e^{x},\,|u_{t}^{N}|^{2}\rangle - i\eta(t)\langle e^{x},\,\overline{u}_{tt}^{N}u_{t}^{N}\rangle
	-i\chi\langle \overline{u}_{xxxt}^{N},\,e^xu_{t}^{N}\rangle
	+ \gamma\mu'(t) \langle \overline{u}_{xx}^{N},\,e^x u_{t}^{N} \rangle\\
	&+ \gamma \mu(t) \langle \overline{u}_{xxt}^{N},\,e^x u_{t}^{N}\rangle+ \overline{\langle A_{t}u_{x}^{N},\,e^{x}\overline{u}_{t}^{N}\rangle} +\overline{ \langle Au_{xt}^{N},\,e^x\overline{u}_{t}^{N}\rangle}
	- i\varepsilon\langle \overline{u}_{xxxxt}^{N},\,e^{x}u_{t}^{N}\rangle  \\
	&-\eta'(t)\overline{\langle |u^{N}|^{2}u^{N},\,e^x\overline{u}_{t}^{N}\rangle}- \eta(t)\overline{\langle (|u^{N}|^{2})_{t}u^{N},\,e^x\overline{u}_{t}^{N}\rangle}
	- \eta(t)\langle e^x,\,|u^{N}|^{4}\rangle=0.
\end{align*}
Subtracting the above equations we have that the terms above are estimate in the same way as previously.
\begin{align}\label{eq14}
	&2i\eta'(t)\langle e^{x},\,|u_{t}^{N}|^{2}\rangle
	+i\eta(t)\left(\langle e^{x},\,\overline{u}_{t}^{N}u_{tt}^{N}\rangle+i\eta(t)\langle e^{x},\,\overline{u}_{tt}^{N}u_{t}^{N}\rangle\right)
	+i\chi\left(\langle u_{xxxt}^{N},\,e^x\overline{u}_{t}^{N}\rangle+\overline{\langle u_{xxxt}^{N},\,e^x\overline{u}_{t}^{N}\rangle}\right)\nonumber\\
	&+\gamma\mu'(t) \left(\langle u_{xx}^{N},\,e^x\overline{u}_{t}^{N} \rangle-\overline{\langle u_{xx}^{N},\,e^x\overline{u}_{t}^{N} \rangle}\right)
	+\gamma\mu(t) \left(\langle u_{xxt}^{N},\,e^x\overline{u}_{t}^{N} \rangle-\overline{\langle u_{xxt}^{N},\,e^x\overline{u}_{t}^{N} \rangle}\right)\nonumber\\
	&+\langle A_{t}u_{x}^{N},\,e^{x}\overline{u}_{t}^{N}\rangle-\overline{\langle A_{t}u_{x}^{N},\,e^{x}\overline{u}_{t}^{N}\rangle}
	+\langle Au_{xt}^{N},\,e^x\overline{u}_{t}^{N}\rangle-\overline{\langle Au_{xt}^{N},\,e^x\overline{u}_{t}^{N}\rangle}\nonumber\\
	&+i\varepsilon\left(\langle u_{xxxxt}^{N},\,e^{x}\overline{u}_{t}^{N}\rangle+\overline{\langle u_{xxxxt}^{N},\,e^{x}\overline{u}_{t}^{N}\rangle}\right)
	-\eta'(t)\left(\langle |u^{N}|^{2}u^{N},\,e^x\overline{u}_{t}^{N}\rangle-\overline{\langle |u^{N}|^{2}u^{N},\,e^x\overline{u}_{t}^{N}\rangle}\right)\nonumber\\
	&- \eta(t)\left(\langle (|u^{N}|^{2})_{t}u^{N},\,e^x\overline{u}_{t}^{N}\rangle-\overline{\langle (|u^{N}|^{2})_{t}u^{N},\,e^x\overline{u}_{t}^{N}\rangle}\right)
	=0.
\end{align}
Note that
\begin{align*}
	&\langle u_{xxxxt}^{N},\,e^{x}\overline{u}_{t}^{N}\rangle+\overline{\langle u_{xxxxt}^{N},\,e^{x}\overline{u}_{t}^{N}\rangle}=
	-\int_0^1u_{xxxt}^N e^x\overline{u}_t^N\ dx-\int_0^1u_{xxxt}^N e^x\overline{u}_{xt}^N\ dx\\
	&\ \ -\int_0^1\overline{u}_{xxxt}^N e^xu_t^N\ dx-\int_0^1\overline{u}_{xxxt}^N e^xu_{xt}^N\ dx
	=\int_0^1u_{xxt}^N e^x\overline{u}_t^N\ dx+\int_0^1e^x\frac{d}{dx}|u_{xt}^N|^2\ dx\\
	& \ \ + \int_0^1e^x\frac{d}{dx}|u_{xt}^N|^2\ dx+2\int_0^1e^x|u_{xxt}^N|^2\ dx+\int_0^1\overline{u}_{xxt}^Ne^x u_{xt}^N\ dx
	=-\int_0^1u_{xt}^N e^x\overline{u}_t^N\ dx-\int_0^1e^x|u_{xt}^N|^2\ dx\\
	& \ \ -\int_0^1e^x |u_{xt}^N|^2\ dx-\int_0^1e^x |u_{xt}^N|^2\ dx +2\int_0^1e^x |u_{xxt}^N|^2\ dx-\int_0^1\overline{u}_{xt}^Ne^xu_t^N\ dx- \int_0^1e^x|u_{xt}^N|^2\ dx\\
	&=-\int_0^1e^x\frac{d}{dx}|u_t^N|^2\ dx-4\langle e^x,\,|u_{xt}^N|^2 \rangle.
\end{align*}
\begin{align*}
	&\langle u_{xxxt}^{N},\,e^x\overline{u}_{t}^{N}\rangle+\overline{\langle u_{xxxt}^{N},\,e^x\overline{u}_{t}^{N}\rangle}
	=-\int_0^1u_{xxt}^N\overline{u}_t^Ne^x\ dx-\int_0^1u_{xxt}^Ne^x\overline{u}_{xt}^N\ dx
	-\int_0^1\overline{u}_{xxt}^N u_t^Ne^x\ dx\\
	&\ \  -\int_0^1\overline{u}_{xxt}^Ne^x\overline{u}_{xt}^N\ dx
	=2\int_0^1e^x|u_{xt}^N|^2\ dx+\int_0^1e^x\frac{d}{dx}|u_{t}^N|^2\ dx-\int_0^1e^x\frac{d}{dx}|u_{xt}^N|^2\ dx\\
	&=3\langle e^x,\,|u_{xt}^{N}|^2\rangle-\langle e^x,\,|u_{t}^{N}|^2\rangle.
\end{align*}
\begin{align*}
	&\langle u_{xx}^{N},\,e^x\overline{u}_{t}^{N} \rangle-\overline{\langle u_{xx}^{N},\,e^x\overline{u}_{t}^{N} \rangle}
	=-\int_0^1e^xu_{x}^N\overline{u}_t^N\ dx-\int_0^1e^xu_{x}^N\overline{u}_{xt}^N\ dx+\int_0^1e^x\overline{u}_{x}^N u_t^N\ dx\\
	&\ \ +\int_0^1e^x\overline{u}_{x}^Nu_{xt}^N\ dx.
\end{align*}
\begin{align*}
	&\langle u_{xxt}^{N},\,e^x\overline{u}_{t}^{N} \rangle-\overline{\langle u_{xxt}^{N},\,e^x\overline{u}_{t}^{N} \rangle}
	=-\int_0^1u_{xt}^Ne^x\overline{u}_t^N\ dx+\int_0^1\overline{u}_{xt}^Ne^xu_t^N\ dx.
\end{align*}
Replacing in \eqref{eq14} we have
\begin{align*}
	&2i\eta'(t)\langle e^{x},\,|u_{t}^{N}|^{2}\rangle+i\eta(t)\frac{d}{dt}\langle e^{x},\,|u_{t}^{N}|^2\rangle
	+3i\chi\langle e^{x},\,|u_{xt}^{N}|^2\rangle-i\chi\langle e^{x},\,|u_{t}^{N}|^2\rangle\\
	&+2i\gamma\mu'(t)\  Im\langle u_{xx}^{N},\,e^x\overline{u}_{t}^{N} \rangle
	+2i\gamma\mu(t)\ Im\langle e^x,\,\overline{u}_{xt}^Nu_{t}^{N} \rangle
	+2i\ Im\langle A_{t}u_{x}^{N},\,e^{x}\overline{u}_{t}^{N}\rangle\\
	&+2i \ Im\langle Au_{xt}^{N},\,e^x\overline{u}_{t}^{N}\rangle+i\varepsilon \langle e^x,\,|u_{t}^{N}|^2\rangle-4i\varepsilon\langle e^x,\,|u_{xt}^{N}|^2\rangle\\
	&-2i\eta'(t)\ Im\langle |u^{N}|^2u^N,\,e^{x}\overline{u}_{t}^{N}\rangle
	-2i\eta(t)\ Im\langle (|u^{N}|^2)_tu^N,\,e^{x}\overline{u}_{t}^{N}\rangle=0.
\end{align*}
Then
\begin{align}\label{eq15}
	&\frac{d}{dt}\left(\eta(t)\langle e^{x},\,|u_{t}^{N}|^{2}\rangle\right)
	+(3\chi-4\varepsilon)\langle e^{x},\,|u_{xt}^{N}|^2\rangle
	+(\varepsilon-\chi)\langle e^{x},\,|u_{t}^{N}|^2\rangle\nonumber\\
	& \leq \ \ -2\gamma\mu'(t)\  Im\langle u_{xx}^{N},\,e^x\overline{u}_{t}^{N} \rangle
	-2\gamma\mu(t)\ Im\langle e^x,\,\overline{u}_{xt}^Nu_{t}^{N} \rangle
	-2\ Im\langle A_{t}u_{x}^{N},\,e^{x}\overline{u}_{t}^{N}\rangle\nonumber\\
	&-2 \ Im\langle Au_{xt}^{N},\,e^x\overline{u}_{t}^{N}\rangle
	+2\eta'(t)\ Im\langle |u^{N}|^2u^N,\,e^{x}\overline{u}_{t}^{N}\rangle\nonumber\\
	&+2\eta(t)\ Im\langle (|u^{N}|^2)_tu^N,\,e^{x}\overline{u}_{t}^{N}\rangle
\end{align}
We are going to bounded the right side of \eqref{eq15}
\begin{align*}
	&-2\gamma\mu'(t)\  Im\langle u_{xx}^N,\,e^{x}\overline{u}_{t}^{N}\rangle\leq
	\varepsilon_1 \|u_{xx}^N\|_{L^2(0,\,1)}^2+C(\varepsilon_1)e\langle e^x,\,|u_{t}^{N}|^2\rangle.\\
	&-2\gamma\mu(t)\  Im\langle e^x,\,\overline{u}_{xt}^{N}u_t^N\rangle\leq
	\varepsilon_2e \langle e^x,\,|\overline{u}_{xt}^{N}|^2\rangle+C(\varepsilon_2)e\langle e^x,\,|u_{t}^{N}|^2\rangle.\\
	&-2\  Im\langle A_{t}u_{x}^{N},\,e^{x}\overline{u}_{t}^{N}\rangle\leq 2\sup_{x\in(0,1)}|A_t|\ \int_0^1|u_{x}^Ne^x\overline{u}_t^N|\ dx
	\leq \varepsilon_3\|u_{x}^N\|_{L^2(0,\,1)}^2+C(\varepsilon_3)e\langle e^x,\,|u_{t}^{N}|^2\rangle.\\
	&-2\  Im\langle A u_{x}^{N},\,e^{x}\overline{u}_{t}^{N}\rangle\leq 2\sup_{x\in(0,1)}|A|\ \int_0^1|u_{xt}^Ne^{2x}\overline{u}_t^N|\ dx
	\leq \varepsilon_4 e\langle e^x,\,|u_{xt}^{N}|^2\rangle+C(\varepsilon_4)e\langle e^x,\,|u_{t}^{N}|^2\rangle.\\
	&2\eta'(t)\  Im\langle |u^{N}|^2u^N,\,e^{x}\overline{u}_{t}^{N}\rangle\leq C\sup_{x\in(0,1)}|u^N|^2\ \int_0^1|u^Ne^x\overline{u}_t^N|\ dx
	\leq \varepsilon_5\|u^N\|_{L^2(0,\,1)}^2\\
	&\ \ +C(\varepsilon_5)e\langle e^x,\,|u_{t}^{N}|^2\rangle.\\
	&2\eta(t)\  Im\langle (|u^{N}|^2)_tu^N,\,e^{x}\overline{u}_{t}^{N}\rangle\leq C\sup_{x\in(0,1)}(|u^N|^2)_t\ \int_0^1|u^Ne^{x}\overline{u}_t^N|\ dx
	\leq \varepsilon_6 \|u^N\|_{L^2(0,\,1)}^2\\
	&\ \ +C(\varepsilon_6)e\langle e^x,\,|u_{t}^{N}|^2\rangle.
\end{align*}
How $3\chi-4\varepsilon>0$ with $\varepsilon_2$ and $\varepsilon_4$ suitable.
Replacing in \eqref{eq15} we have
\begin{align}\label{eq17}
	&\frac{d}{dt}\left(\eta(t)\langle e^{x},\,|u_{t}^{N}|^{2}\rangle\right)
	+(3\chi-4\varepsilon-e\varepsilon_2-e\varepsilon_4)\langle e^{x},\,|u_{xt}^{N}|^2\rangle
	\leq \ \ C\langle e^x,\,|u_{t}^{N}|^2\rangle\nonumber \\
	& \ \ +\varepsilon_1\|u_{xx}^N\|_{L^2(0,\,1)}^2+\varepsilon_3\|u_{x}^N\|_{L^2(0,\,1)}^2
	+C\|u^N\|_{L^2(0,\,1)}^2.
\end{align}
Integrating over $(0,t)$ with $t\leq T$.
\begin{align*}
	&\eta(t)\langle e^{x},\,|u_{t}^{N}|^{2}\rangle \ \ \leq \ \  \eta(0)\langle e^{x},\,|u_{t}^{N}(\cdot,\,0)|^{2}\rangle
	+\varepsilon_1\int_0^T\|u_{xx}^N\|_{L^2(0,\,1)}^2ds\\
	&+\varepsilon_3\int_0^T\|u_{x}^N\|_{L^2(0,\,1)}^2ds
	+C\int_0^T\|u^N\|_{L^2(0,\,1)}^2ds+C\int_0^t\langle e^{x},\,|u_{t}^{N}|^{2}\rangle ds.
\end{align*}
On the other hand  of \eqref{307} with $t=0$ and $g_j=u_t^N(\cdot,\,0)$ we have
\begin{align*}
	&i\|u_t^N(\cdot,\,0)\|_{L^2(0,\,1)}^2 + i\chi\mathfrak{B}(0)\langle u_{xxx}^{N}(\cdot,\,0),\,\overline{u}_t^N(\cdot,\,0)\rangle + \gamma \mathfrak{C}(0)\langle u_{xx}^{N}(\cdot,\,0),\,\overline{u}_t^N(\cdot,\,0)\rangle \\
	&+\langle I(\cdot,\,0)\ u_{x}^{N}(\cdot,\,0),\,\overline{u}_t^N(\cdot,\,0)\rangle
	+ i\varepsilon\mathfrak{B}(0)\langle u_{xxxx}^{N}(\cdot,\,0),\,\overline{u}_t^N(\cdot,\,0)\rangle
	- \langle |u^{N}(\cdot,\,0)|^{2}u^{N}(\cdot,\,0),\,\overline{u}_t^N(\cdot,\,0)\rangle\\
	& =  0
\end{align*}
then
\begin{align*}
	&\|u_t^N(\cdot,\,0)\|_{L^2(0,\,1)}^2\leq   \varepsilon_8\|u_{xxx}^N(\cdot,\,0)\|_{L^2(0,\,1)}^2+C(\varepsilon_8)\|u_t^N(\cdot,\,0)\|_{L^2(0,\,1)}^2\\
	& \ \ +\varepsilon_9\sup_{x\in(0,1)}|I(\cdot,\,0)|\|u_{x}^N(\cdot,\,0)\|_{L^2(0,\,1)}^2 +\sup_{x\in(0,1)}|I(\cdot,\,0)| \ C(\varepsilon_9)\|u_t^N(\cdot,\,0)\|_{L^2(0,\,1)}^2\\
	& \ \ + \varepsilon_{10}\|u_{xxxx}^N(\cdot,\,0)\|_{L^2(0,\,1)}^2+C(\varepsilon_{10})\|u_t^N(\cdot,\,0)\|_{L^2(0,\,1)}^2\\
	& \ \ +\sup_{x\in(0,1)}|u_0^N|^2 \ \varepsilon_{11}\|u^N(\cdot,\,0)\|_{L^2(0,\,1)}^2 +\sup_{x\in(0,1)}|u_0^N|^2\ C(\varepsilon_{11})\|u_t^N(\cdot,\,0)\|_{L^2(0,\,1)}^2
\end{align*}
then for $C(\varepsilon_i)$ with $i=8,\ldots,11$ suitable, we have
\begin{align}\label{eq30}
	\|u_t^N(\cdot,\,0)\|_{L^2(0,\,1)}^2\leq& C\left(\|u_0^N\|_{H^4(0,\,1)}+\|J(\cdot,\,0)\|_{L^2(0,\,1)}^2\right)
\end{align}
\begin{remark}
	If $\varepsilon=0$ then \eqref{eq30} is true, therefore the constant C in \eqref{eq30} independent of $\varepsilon$.
\end{remark}
Then
\begin{align}\label{eq28}
	&\langle e^{x},\,|u_{t}^{N}|^{2}\rangle \ \ \leq \ \  C+C\int_0^t\langle e^{x},\,|u_{t}^{N}|^{2}\rangle ds.
\end{align}
By Gronwall inequality
\begin{align}\label{eq18}
	\|u_t^N\|_{L^2(0,\,1)}^2\leq \langle e^{x},\,|u_{t}^{N}|^{2}\rangle\leq C\exp(Ct)
\end{align}
for all $t\in (0,T)$ where the constant $C$ independent of $\varepsilon$. Therefore
\begin{align}\label{eq29}
	(u_t^N) \ \ \textrm{is bounded in} \ \  L^\infty(0,T:L^2(0,1))\hookrightarrow L^2(0,T:L^2(0,1))
\end{align}
independent of $\varepsilon$. Note that of \eqref{eq17} together with  \eqref{eq18} we have
\begin{align}\label{eq19}
	(3\chi-4\varepsilon-e\varepsilon_2-e\varepsilon_4)\int_0^t\langle e^{x},\,|u_{xt}^{N}|^2\rangle\ ds
	\leq  C.
\end{align}
\begin{remark}
	If $\varepsilon=0$, then for $\varepsilon_2$ and $\varepsilon_4$ suitable, we have \eqref{eq19} is true, where $C$ independent of $\varepsilon$.
\end{remark}
Therefore
\begin{align}\label{eq31}
	(u_{xt}^N) \ \ \textrm{is bounded in} \ \   L^2(0,T:L^2(0,1)) \ \ \textrm{indepent of }\varepsilon
\end{align}
{\it Step four.} We obtain a bound for $u_{xxt}^N$.  Differenting \eqref{407} respect $t$, and replacing $g_{j}$ by $\overline{u}_{t}^{N}$, we have
\begin{align*}
	&i\eta(t)\langle u_{tt}^{N},\,\overline{u}_{t}^{N}\rangle + i\eta'(t)\|u_{t}^{N}(t)\|_{L^{2}(0,\,1)}^{2} + i\chi\langle u_{xxxt}^{N},\,\overline{u}_{t}^{N}\rangle  \nonumber \\
	& + \gamma \mu(t)\langle u_{xxt}^{N},\,\overline{u}_{t}^{N}\rangle+\gamma \mu'(t)\langle u_{xx}^{N},\,\overline{u}_{t}^{N}\rangle
	+ \langle A_{t}u_{x}^{N},\,\overline{u}_{t}^{N}\rangle + \langle Au_{xt}^{N},\,\overline{u}_{t}^{N}\rangle  \nonumber \\
	&+ i\varepsilon\langle u_{xxxxt}^{N},\,\overline{u}_{t}^{N}\rangle
	- \eta(t)\langle (|u^{N}|^{2})_{t}u^{N},\,\overline{u}_{t}^{N}\rangle
	- \eta(t)\langle |u^{N}|^{2}u_{t}^{N},\, \overline{u}_{t}^{N}\rangle \nonumber \\
	&- \eta'(t)\langle |u^{N}|^{2}u^{N},\,\overline{u}_{t}^{N}\rangle =  0.
\end{align*}
Hence
\begin{align*}
	&i\eta(t)\langle u_{tt}^{N},\,\overline{u}_{t}^{N}\rangle + i\eta'(t)\|u_{t}^{N}(t)\|_{L^{2}(0,\,1)}^{2} - i\chi\langle u_{xxt}^{N},\,\overline{u}_{xt}^{N}\rangle  \nonumber \\
	& - \gamma\mu(t)\|u_{xt}^{N}(t)\|_{L^{2}(0,\,1)}^{2} - \gamma \mu'(t)\langle u_{x}^{N},\,\overline{u}_{xt}^{N}\rangle + \langle A_{t}u_{x}^{N},\,\overline{u}_{t}^{N}\rangle + \langle Au_{xt}^{N},\,\overline{u}_{t}^{N}\rangle \nonumber \\
	&+ i\varepsilon\langle u_{xxxxt}^{N},\,\overline{u}_{t}^{N}\rangle
	- \eta(t)\langle (|u^{N}|^{2})_{t}u^{N},\,\overline{u}_{t}^{N}\rangle
	- \eta(t)\langle |u^{N}|^{2},\,|u_{t}^{N}|^{2}\rangle - \eta'(t)\langle |u^{N}|^{2}u^{N},\,\overline{u}_{t}^{N}\rangle =  0.
\end{align*}
Applying conjugate
\begin{align*}
	&-i\eta(t)\langle \overline{u}_{tt}^{N},\,u_{t}^{N}\rangle - i\eta'(t)\|u_{t}^{N}(t)\|_{L^{2}(0,\,1)}^{2} + i\chi\langle \overline{u}_{xxt}^{N},\,u_{xt}^{N}\rangle  \nonumber \\
	& - \gamma \mu(t)\|u_{xt}^{N}(t)\|_{L^{2}(0,\,1)}^{2} - \gamma \mu'(t)\langle \overline{u}_{x}^{N},\,u_{xt}^{N}\rangle + \overline{\langle A_{t}u_{x}^{N},\,\overline{u}_{t}^{N}\rangle }+\overline{\langle Au_{xt}^{N},\,\overline{u}_{t}^{N}\rangle} \nonumber \\
	&- i\varepsilon\overline{\langle u_{xxxxt}^{N},\,\overline{u}_{t}^{N}\rangle}
	- \eta(t)\overline{\langle (|u^{N}|^{2})_{t}u^{N},\,\overline{u}_{t}^{N}\rangle}
	- \eta(t)\langle |u^{N}|^{2},\,|u_{t}^{N}|^{2}\rangle - \eta'(t)\overline{\langle |u^{N}|^{2}u^{N},\,\overline{u}_{t}^{N}\rangle} = 0
\end{align*}
Subtracting the above equations yield
\begin{align}\label{eq5}
	&i\eta(t)\frac{d}{dt}\|u_{t}^{N}\|_{L^2(0,\,1)}^2 +2i\eta'(t)\|u_{t}^{N}\|_{L^{2}(0,\,1)}^{2}
	- i\chi\left(\langle u_{xxt}^{N},\, \overline{u}_{xt}^{N}\rangle +\langle \overline{u}_{xxt}^{N},\,u_{xt}^{N}\rangle\right) \nonumber\\
	& - \gamma\mu'(t)\left(\langle u_{x}^{N},\,\overline{u}_{xt}^{N}\rangle-\langle \overline{u}_{x}^{N},\,u_{xt}^{N}\rangle\right)
	+2iIm\langle A_{t}u_{x}^{N},\,\overline{u}_{t}^{N}\rangle +2iIm\langle Au_{xt}^{N},\,\overline{u}_{t}^{N}\rangle \nonumber \\
	&+i\varepsilon\left(\langle u_{xxxxt}^{N},\,\overline{u}_{t}^{N}\rangle+\langle \overline{u}_{xxxxt}^{N},\,u_{t}^{N}\rangle\right)
	- \eta(t)\left(\langle (|u^{N}|^{2})_{t}u^{N},\,\overline{u}_{t}^{N}\rangle-\overline{\langle (|u^{N}|^{2})_{t}u^{N},\,\overline{u}_{t}^{N}\rangle}\right)\nonumber\\
	&- \eta'(t)\left(\langle |u^{N}|^{2}u^{N},\,\overline{u}_{t}^{N}\rangle-\overline{\langle |u^{N}|^{2}u^{N},\,\overline{u}_{t}^{N}\rangle}\right) = 0.
\end{align}
Note that
\begin{align*}
	\langle u_{xxxxt}^{N},\,\overline{u}_{t}^{N}\rangle+\langle \overline{u}_{xxxxt}^{N},\,u_{t}^{N}\rangle=2\|u_{xxt}^N\|_{L^2(0,\,1)}^2
\end{align*}
then replacing in \eqref{eq5} we have
\begin{align*}
	&i\eta(t)\frac{d}{dt}\|u_{t}^{N}\|_{L^{2}(0,\,1)}^{2} + 2i\eta'(t)\|u_{t}^{N}\|_{L^{2}(0,\,1)}^{2} - i\chi\int_0^1\frac{d}{dx}|u_{xt}|^2dx\nonumber \\
	&-2i\gamma\mu'(t)Im\langle u_{x}^{N},\,\overline{u}_{xt}^{N}\rangle + 2iIm\langle A_{t}u_{x}^{N},\,\overline{u}_{t}^{N}\rangle + 2iIm\langle Au_{xt}^{N},\,\overline{u}_{t}^{N}\rangle \nonumber \\
	&+2i\chi|u_{xt}(1,\,t)|^2+2i\varepsilon\|u_{xxt}^{N}\|_{L^{2}(0,\,1)}^{2}-2i\eta(t) Im\langle (|u^{N}|^{2})_tu^{N},\,\overline{u}_{t}^{N}\rangle\\
	& - 2i\eta'(t)Im\langle |u^{N}|^{2}u^{N},\,\overline{u}_{t}^{N}\rangle =  0
\end{align*}
or
\begin{align}\label{eq6}
	&\eta(t)\frac{d}{dt}\|u_{t}^{N}\|_{L^{2}(0,\,1)}^{2} + 2\eta'(t)\|u_{t}^{N}\|_{L^{2}(0,\,1)}^{2}
	+\chi|u_{xt}(1,\,t)|^2+\chi|u_{xt}(0,\,t)|^2+2\varepsilon\|u_{xxt}^{N}\|_{L^{2}(0,\,1)}^{2}\nonumber\\
	=&2\gamma\mu'(t)Im\langle u_{x}^{N},\,\overline{u}_{xt}^{N}\rangle
	-2Im\langle A_{t}u_{x}^{N},\,\overline{u}_{t}^{N}\rangle - 2Im\langle Au_{xt}^{N},\,\overline{u}_{t}^{N}\rangle\nonumber\\
	&+2\eta(t) Im\langle (|u^{N}|^{2})_tu^{N},\,\overline{u}_{t}^{N}\rangle
	+2\eta'(t)Im\langle |u^{N}|^{2}u^{N},\,\overline{u}_{t}^{N}\rangle.
\end{align}
We are going to bounded the right side of \eqref{eq6}
\begin{align}\label{eq7}
	&2\gamma\mu'(t)Im\langle u_{x}^{N},\,\overline{u}_{xt}^{N}\rangle\leq C\left(\|u_{x}^{N}\|_{L^{2}(0,\,1)}^{2}+\|u_{xt}^{N}\|_{L^{2}(0,\,1)}^{2}\right)\nonumber\\
	&-2Im\langle A_{t}u_{x}^{N},\,\overline{u}_{t}^{N}\rangle \leq \sup_{x\in(0,1)}|A_t|\|u_{x}^{N}\|_{L^{2}(0,\,1)}^{2}+\sup_{x\in(0,1)}|A_t|\|u_{t}^{N}\|_{L^{2}(0,\,1)}^{2}\nonumber\\
	&- 2Im\langle Au_{xt}^{N},\,\overline{u}_{t}^{N}\rangle\leq
	\sup_{x\in(0,1)}|A|\|u_{xt}^{N}\|_{L^{2}(0,\,1)}^{2}+\sup_{x\in(0,1)}|A|\|u_{t}^{N}\|_{L^{2}(0,\,1)}^{2}\nonumber\\
	&2\eta(t) Im\langle (|u^{N}|^{2})_tu^{N},\,\overline{u}_{t}^{N}\rangle\leq
	C\left(\sup_{x\in(0,1)}|(|u^N|^2)_t|\ \|u^{N}\|_{L^{2}(0,\,1)}^{2}+\sup_{x\in(0,1)}|(|u^N|^2)_t|\ \|u_t^{N}\|_{L^{2}(0,\,1)}^{2}\right)\nonumber\\
	&2\eta'(t)Im\langle |u^{N}|^{2}u^{N},\,\overline{u}_{t}^{N}\rangle\leq
	C\left(\sup_{x\in(0,1)}|u^N|^2\ \|u^{N}\|_{L^{2}(0,\,1)}^{2}+\sup_{x\in(0,1)}|u^N|^2\ \|u_t^{N}\|_{L^{2}(0,\,1)}^{2}\right).
\end{align}
Integrating \eqref{eq6} over 0 to $t$ with $t\leq T$ together with \eqref{eq7}  we have
\begin{align}\label{eq8}
	&\eta(t)\|u_{t}^{N}\|_{L^{2}(0,\,1)}^{2}
	+2\varepsilon\int_0^t\|u_{xxt}^{N}\|_{L^{2}(0,\,1)}^{2}ds\leq\eta(0) \|u_{t}^{N}(\cdot,\,0)\|_{L^{2}(0,\,1)}^{2}
	+C\int_0^t\|u_{t}^{N}\|_{L^{2}(0,\,1)}^{2}ds\nonumber\\
	&\ \ +C\int_0^t\|u_{x}^{N}\|_{L^{2}(0,\,1)}^{2}ds+C\int_0^t\|u_{xt}^{N}\|_{L^{2}(0,\,1)}^{2}ds+C\int_0^t\|u^{N}\|_{L^{2}(0,\,1)}^{2}ds.
\end{align}
Since \eqref{eq8} we obtain
\begin{align}\label{eq9}
	&\eta(t)\|u_{t}^{N}\|_{L^{2}(0,\,1)}^{2}
	+2\varepsilon\int_0^T\|u_{xxt}^{N}\|_{L^{2}(0,\,1)}^{2}ds\nonumber\\
	&\leq C\|u_{0}^{N}\|_{H^4(0,\,1)}^2
\end{align}
therefore
\begin{align}
	\label{eq56}(u_{xxt}^N) \ \ \textrm{is bounded in} \ \  L^2(0,T:L^2(0,1)) \ \ \textrm{dependent of }\varepsilon\\
\end{align}
{\it Step to the limit.}
Now follows  summarize our results.
\begin{align}
	&\label{eq32}(u^N) \ \ \textrm{is bounded in} \ \   L^\infty(0,T:L^2(0,1))\hookrightarrow L^2(0,T:L^2(0,1)) \ \ \textrm{independent of }\varepsilon\\
	&\label{eq33}(u_x^N) \ \ \textrm{is bounded in} \ \  L^2(0,T:L^2(0,1)) \ \ \textrm{indepent of }\varepsilon\\
	&\label{eq34}(u_{xx}^N) \ \ \textrm{is bounded in} \ \  L^\infty(0,T:L^2(0,1))\hookrightarrow L^2(0,T:L^2(0,1)) \ \ \textrm{independent of }\varepsilon\\
	&\label{eq35}(u_{xxxx}^N) \ \ \textrm{is bounded in} \ \  L^2(0,T:L^2(0,1)) \ \ \textrm{independent of }\varepsilon\\
	&\label{eq36}(u_t^N) \ \ \textrm{is bounded in} \ \  L^\infty(0,T:L^2(0,1))\hookrightarrow  L^2(0,T:L^2(0,1)) \ \ \textrm{independent of }\varepsilon\\
	&\label{eq37}(u_{xt}^N) \ \ \textrm{is bounded in} \ \ L^2(0,T:L^2(0,1)) \ \ \textrm{independent of }\varepsilon\\
	&\label{eq38}(u_{xxt}^N) \ \ \textrm{is bounded in} \ \   L^2(0,T:L^2(0,1)) \ \ \textrm{dependent of }\varepsilon
\end{align}
We affirm that
\begin{align}\label{eq39}
	(u_x^N) \ \  \textrm{is bounded in} \ \  L^\infty(0,T:H^1(0,1))\ \   \textrm{independent of} \ \  \varepsilon.
\end{align}
In effect, by the Gagliardo-Nirenberg inequality, we have
$$\|u_x^N\|_{L^2(0,1)}\leq\alpha\|u_{xx}^N\|_{L^2(0,1)}+C(\alpha)\|u^N\|_{L^2(0,1)}$$
for all $t\in[0,T]$ then
$$\sup_{t\in[0,T]}\|u_x^N\|_{L^2(0,1)}\leq\alpha\sup_{t\in[0,T]}\|u_{xx}^N\|_{L^2(0,1)}+C(\alpha)\sup_{t\in[0,T]}\|u^N\|_{L^2(0,1)}<C$$
the result follows.\\
Gathering the estimates \eqref{eq32}-\eqref{eq38} we have there exists an subsequences of $(u^N)$ and $(u_t^N)$ that also denote by $(u^N)$ and $u_t^N$ such that
\begin{align}
&\label{eq61}u^N\rightharpoonup u, \ \ \textrm{weak}, \ \ u\in L^2(0,T:L^2(0,1))\\
&\label{eq62}u^N\rightharpoonup^* u, \ \ \textrm{star weak}, \ \ u\in L^\infty(0,T:L^2(0,1))\\
&\label{eq63}u_x^N\rightharpoonup u_x, \ \ \textrm{weakly}, \ \ u_x\in L^2(0,T:L^2(0,1))\\
&\label{eq64}u_x^N\rightharpoonup^* u_x, \ \ \textrm{star weak}, \ \ u_x\in L^\infty(0,T:L^2(0,1))\\
&\label{eq65}u_t^N\rightharpoonup u_t, \ \ \textrm{weak}, \ \ u_t\in L^2(0,T:L^2(0,1))\\
&\label{eq66}u_t^N\rightharpoonup^* u_t, \ \ \textrm{star weak}, \ \ u_t\in L^\infty(0,T:L^2(0,1))
\end{align}
the convergence independent of $\varepsilon$.\\
Finally multiplying  \eqref{301} by $\varphi\in C_0^\infty(0,1)$ test functions and integrating over $(0,1)$, we have
\begin{align*}
	&i\int_0^1u_{t}\varphi \ dx + i\chi\mathfrak{B}(t)\int_0^1u_{xxx}\varphi \ dx
	+ \gamma  \mathfrak{C}(t)\int_0^1u_{xx}\varphi \ dx
	+\int_0^1 I(x,\,t)u_{x}\varphi \ dx \\
	&+ i\varepsilon \mathfrak{B}(t)\int_0^1u_{xxxx}\varphi \ dx - \int_0^1|u|^{2}u\varphi \ dx = 0.
\end{align*}
Integrating of parts, we have
\begin{align}\label{eq40}
	&i\langle u_{t},\,\varphi\rangle + i\chi\mathfrak{B}(t)\langle u_x,\,\frac{d^2\varphi}{dx^2}\rangle
	- \gamma  \mathfrak{C}(t)\langle u_{x},\,\frac{d\varphi}{dx}\rangle
	+\langle I(x,\,t)u_{x},\,\varphi \rangle \\
	&- i\varepsilon \mathfrak{B}(t)\langle u_{x},\,\frac{d^3\varphi}{dx^3} \rangle - \langle |u|^{2}u,\,\varphi\rangle = 0.
\end{align}
for all $\varphi\in C_0^\infty(0,1)$.\\
Since $u^N$ es approximation solution, we have
\begin{align}\label{eq41}
	&i\langle u_{t}^N,\,\varphi\rangle + i\chi\mathfrak{B}(t)\langle u_x^N,\,\frac{d^2\varphi}{dx^2}\rangle
	- \gamma  \mathfrak{C}(t)\langle u_{x}^N,\,\frac{d\varphi}{dx}\rangle
	+\langle I(x,\,t)u_{x}^N,\,\varphi \rangle \nonumber \\
	&- i\varepsilon \mathfrak{B}(t)\langle u_{x}^N,\,\frac{d^3\varphi}{dx^3} \rangle - \langle |u^N|^{2}u^N,\,\varphi\rangle = 0.
\end{align}
Since $H^1(0,1)\subset C([0,1])\subset L^2(0,1)$ with the first injection  compact and the second injection continuous, together with \eqref{eq32}-\eqref{eq39}, by the Aubin -Lions lemma and Simon's compactness result, we have there exists
$u\in C([0,T]:C([0,1]))$ and $u_x\in C([0,T]:C([0,1]))$ such that a subsequence of $(u^N)$  and $(u_x^N)$ that denote also by $(u^N)$ and $(u_x^N)$ that independent of $\varepsilon$, that converge (strong) to $u$ and $u_x$ en $C([0,T]:C([0,1]))$. This is
\begin{align}\label{eq42}
	\|u^N-u\|_{C([0,T]:C([0,1]))}\rightarrow 0\hspace{1cm}\|u_x^N-u_x\|_{C([0,T]:C([0,1]))}\rightarrow 0.
\end{align}
It follows  that approximations $u^{N}(x,\,t)$ can be extended to entire $Q = (0,\,1)\times (0,\,T)$ with arbitrary finite $T > 0$, gathering with \eqref{eq61}-\eqref{eq66} and justify the passage to the limit as $N\rightarrow +\infty$ in \eqref{eq41}. This is
\begin{align}\label{eq43}
	&i\langle u_{t},\,\varphi\rangle + i\chi\mathfrak{B}(t)\langle u_x,\,\frac{d^2\varphi}{dx^2}\rangle
	- \gamma  \mathfrak{C}(t)\langle u_{x},\,\frac{d\varphi}{dx}\rangle
	+\langle I(x,\,t)u_{x},\,\varphi \rangle \nonumber \\
	&- i\varepsilon \mathfrak{B}(t)\langle u_{x},\,\frac{d^3\varphi}{dx^3} \rangle - \langle |u|^{2}u,\,\varphi\rangle = 0
\end{align}
The solution independent of $\varepsilon$, take $\varepsilon\rightarrow0$, we have
\begin{align}\label{eq44}
	&i\langle u_{t},\,\varphi\rangle + i\chi\mathfrak{B}(t)\langle u_x,\,\frac{d^2\varphi}{dx^2}\rangle
	- \gamma  \mathfrak{C}(t)\langle u_{x},\,\frac{d\varphi}{dx}\rangle
	+\langle I(x,\,t)u_{x},\,\varphi \rangle - \langle |u|^{2}u,\,\varphi\rangle = 0
\end{align}
for all $\varphi\in C_0^\infty(0,1)$.
Since $H^4(0,1)\hookrightarrow C^3([0,1])$ together with \eqref{eq35} we have $(u^N)\subset C^3[0,1]$. Then since \eqref{eq44} we have
\begin{align}\label{eq45}
	&i\langle u_{t},\,\varphi\rangle + i\chi\mathfrak{B}(t)\langle u_{xxx},\,\varphi\rangle
	- \gamma  \mathfrak{C}(t)\langle u_{xx},\,\varphi\rangle
	+\langle I(x,\,t)u_{x},\,\varphi \rangle - \langle |u|^{2}u,\,\varphi\rangle = 0
\end{align}
for all $\varphi\in C_0^\infty(0,1)$, then
\begin{align}\label{eq46}
	&iu_{t} + i\chi\mathfrak{B}(t) u_{xxx}^N
	- \gamma  \mathfrak{C}(t)u_{xx}
	+I(x,\,t)u_{x} - |u|^{2}u = 0.
\end{align}
This prove the existence of a classic solution to \eqref{eq47}-\eqref{eq50}.\\
\renewcommand{\theequation}{\thesection.\arabic{equation}}
\setcounter{equation}{0}
\section{Uniqueness of solutions}
\noindent
The following theorem shows the uniqueness of solutions.
\begin{theorem} There is only one solution of  \eqref{eq47}-\eqref{eq50}.
\end{theorem}
\proof
Assume that $u$ and $w$ are two solutions of \eqref{eq47}-\eqref{eq50}. Therefore
\begin{align}
	\label{501}&i u_{t} + i\chi \mathfrak{B}(t) u_{xxx} + \gamma \mathfrak{C}(t)u_{xx} + I(x,\,t)u_{x} - |u|^{2}u = 0 \\
	\label{502}&i w_{t} + i\chi \mathfrak{B}(t) w_{xxx} + \gamma \mathfrak{C}(t)w_{xx} + I(x,\,t)w_{x} - |w|^{2}w = 0
\end{align}
Subtracting \eqref{501} with \eqref{502} we have
\begin{align}\label{503}
	&i(u - w)_{t} + i\chi\mathfrak{B}(t) (u - w)_{xxx} + \gamma \mathfrak{C}(t)(u - w)_{xx}
	+(u-w)_x I(x,\,t)\nonumber\\
	&- \left[|u|^{2}u- |w|^{2}w\right] = 0.
\end{align}
But
\begin{align*}
	|u|^{2}u - |w|^{2}w = |u|^{2}(u - w) + \left[|u|^{2} - |w|^{2}\right]w.
\end{align*}
Multiplying \eqref{503} by $e^{x}(\overline{u - w})$ and integrating over $(0,\,1)$ and using \eqref{eq48}-\eqref{eq50} it follows that
\begin{align}
	&i\int_{0}^{1}e^{x}(\overline{u - w})(u - w)_{t}dx + i\chi\mathfrak{B}(t)\int_{0}^{1}e^{x}(\overline{u - w})(u - w)_{xxx}dx\nonumber\\
	& + \gamma \mathfrak{C}(t)\int_{0}^{1}e^{x}(\overline{u - w})(u - w)_{xx}dx  +\int_{0}^{1}e^{x}I(x,\,t)(\overline{u - w})(u - w)_{x}dx  \nonumber \\
	\label{505}&-\int_{0}^{1}e^{x}(\overline{u - w})\left(|u|^2(u-w)+(|u|^2 - |w|^2)w\right) dx = 0.
\end{align}
Each term is treated separately integrating by parts and using \eqref{eq47}-\eqref{eq50}. In fact
\begin{align*}
	& i\chi\mathfrak{B}(t)\int_{0}^{1}e^{x}(\overline{u - w})(u - w)_{xxx}dx = i\chi\mathfrak{B}(t)\int_{0}^{1}e^{x}(\overline{u - w})(u - w)_{x}dx \nonumber \\
	&+ i\chi\mathfrak{B}(t)\int_{0}^{1}e^{x}|(u - w)_{x}|^{2}dx - i\chi\mathfrak{B}(t)\int_{0}^{1}e^{x}(\overline{u - w})_{x}(u - w)_{xx}dx
\end{align*}
\begin{align*}
	\gamma \mathfrak{C}(t)\int_{0}^{1}e^{x}(\overline{u - w})(u - w)_{xx}dx
	= -\gamma \mathfrak{C}(t)\int_{0}^{1}e^{x}(\overline{u - w})(u - w)_{x}dx - \gamma \mathfrak{C}(t)\int_{0}^{1}e^{x}|(u - w)_{x}|^{2}dx.
\end{align*}
Replacing into \eqref{505} we have
\begin{align*}
	&i\int_{0}^{1}e^{x}(\overline{u - w})(u - w)_{t}dx + i\chi\mathfrak{B}(t)\int_{0}^{1}e^{x}(\overline{u - w})(u - w)_{x}dx \nonumber \\
	&+ i\chi\mathfrak{B}(t)\int_{0}^{1}e^{x}|(u - w)_{x}|^{2}dx - i\chi\mathfrak{B}(t)\int_{0}^{1}e^{x}(\overline{u - w})_{x}(u - w)_{xx}dx \nonumber \\
	&-\gamma \mathfrak{C}(t)\int_{0}^{1}e^{x}(\overline{u - w})(u - w)_{x}dx
	- \gamma \mathfrak{C}(t)\int_{0}^{1}e^{x}|(u - w)_{x}|^{2}dx \nonumber \\
	& - \int_{0}^{1}e^{x}|u|^{2}|u - w|^{2}dx + \int_{0}^{1}e^{x}I(x,\,t)(\overline{u - w})(u - w)_{x}dx  \nonumber \\
	&- \int_{0}^{1}e^{x}(|u| + |w|)(|u| - |w|)(\overline{u - w})w \ dx = 0.
\end{align*}
Applying conjugate
\begin{align*}
	&-i\int_{0}^{1}e^{x}(u - w)(\overline{u - w})_{t}dx - i\chi\mathfrak{B}(t)\int_{0}^{1}e^{x}(u - w)(\overline{u - w})_{x}dx \nonumber \\
	&- i\chi\mathfrak{B}(t)\int_{0}^{1}e^{x}|(u - w)_{x}|^{2}dx + i\chi\mathfrak{B}(t)\int_{0}^{1}e^{x}(u - v)_{x}(\overline{u - w})_{xx}dx \nonumber \\
	&-\gamma \mathfrak{C}(t)\int_{0}^{1}e^{x}(u - w)(\overline{u - w})_{x}dx
	- \gamma \mathfrak{C}(t)\int_{0}^{1}e^{x}|(u - w)_{x}|^{2}dx \nonumber \\
	& - \int_{0}^{1}e^{x}|u|^{2}|u - w|^{2}dx + \int_{0}^{1}e^{x}I(x,\,t)(u - w)(\overline{u - w})_{x}dx  \nonumber \\
	&- \int_{0}^{1}e^{x}(|u| + |w|)(|u| - |w|)(u - w)\overline{w}dx = 0.
\end{align*}
Then subtracting the above two terms
\begin{align*}
	&\frac{d}{dt}\int_{0}^{1}e^{x}|u - w|^{2}dx + \chi\mathfrak{B}(t)\int_{0}^{1}e^{x}\frac{d}{dx}|u - w|^{2}dx \nonumber \\
	&- \chi\mathfrak{B}(t)\int_{0}^{1}e^{x}\frac{d}{dx}|(u - w)_{x}|^{2}dx + 2\chi\mathfrak{B}(t)\int_{0}^{1}e^{x}|(u - w)_{x}|^{2}dx    \nonumber \\
	&- 2\gamma \mathfrak{C}(t)Im\int_{0}^{1}e^{x}(\overline{u - w})(u - w)_{x}dx + 2Im\int_{0}^{1}e^{x}I(x,\,t)(\overline{u - w})(u - w)_{x}dx  \nonumber \\
	& - 2Im\int_{0}^{1}e^{x}(\overline{u - w})(|u| + |w|)(|u| - |w|)wdx = 0
\end{align*}
Integrating by parts and using \eqref{eq48}-\eqref{eq50} we obtain
\begin{align*}
	&\frac{d}{dt}\int_{0}^{1}e^{x}|u - w|^{2}dx - \chi\mathfrak{B}(t)\int_{0}^{1}e^{x}|u - w|^{2}dx+3\chi\mathfrak{B}(t)\int_{0}^{1}e^{x}|(u - w)_{x}|^{2}dx   \nonumber \\
	&- 2\gamma \mathfrak{C}(t)Im\int_{0}^{1}e^{x}(\overline{u - w})(u - w)_{x}dx + 2Im\int_{0}^{1}e^{x}I(x,\,t)(\overline{u - w})(u - w)_{x}dx  \nonumber \\
	& - 2Im\int_{0}^{1}e^{x}(\overline{u - w})(|u| + |w|)(|u| - |w|)w \ dx = 0.
\end{align*}
or
\begin{align}\label{eq58}
	&\frac{d}{dt}\langle e^{x},\,|u - w|^{2}\rangle - \chi\mathfrak{B}(t)\langle e^{x},\ ,|u - w|^{2}\rangle+3\chi\mathfrak{B}(t)\langle e^{x},\,|(u - w)_{x}|^{2}\rangle   \nonumber \\
	&= 2\gamma \mathfrak{C}(t)Im\int_{0}^{1}e^{x}(\overline{u - w})(u - w)_{x}dx - 2Im\int_{0}^{1}e^{x}I(x,\,t)(\overline{u - w})(u - w)_{x}dx  \nonumber \\
	& + 2Im\int_{0}^{1}e^{x}(\overline{u - w})(|u| + |w|)(|u| - |w|)w \ dx.
\end{align}
We are going to bounded the right side of \eqref{eq58}
\begin{align*}
	2\gamma \mathfrak{C}(t)Im\int_{0}^{1}e^{x}(\overline{u - w})(u - w)_{x}dx &\leq  C\int_{0}^{1}e^{2x}|\overline{u - w}| \ |(u - w)_{x}| \ dx\\
	&\leq \varepsilon_1e\langle e^x,\ ,|(u - w)_{x}|^2\rangle+C(\varepsilon_1)e\langle e^x,\,|u - w|^2\rangle.\\
	- 2Im\int_{0}^{1}e^{x}I(x,\,t)(\overline{u - w})(u - w)_{x}dx&\leq 2\sup_{x\in(0,1)}|I(x,\,t)|\int_{0}^{1}e^{2x}|u - w| \ |(u - w)_{x}| \ dx\\
	&\leq\varepsilon_2e\langle e^x,\ ,|(u - w)_{x}|^2\rangle+C(\varepsilon_2)e\langle e^x,\,|u - w|^2\rangle.\\
	2Im\int_{0}^{1}e^{x}(\overline{u - w})(|u| + |w|)(|u| - |w|)w \ dx&\leq 2(\|u\|_{C(0,T:C([0,1]))}+\|w\|_{C(0,T:C([0,1]))})\\
	&\times \ \ \|w\|_{C(0,T:C([0,1]))}\int_0^1e^x|u-w|^2 \ dx\\
	&\leq C\langle e^x,\, |u-w|^2\rangle.
\end{align*}
We replacing in  \eqref{eq58}, we have
\begin{align}\label{eq59}
	&\frac{d}{dt}\langle e^{x},\,|u - w|^{2}\rangle + (3\chi\mathfrak{B}(t)-\varepsilon_1e-e\varepsilon_2)\langle e^{x},\ ,|(u - w)_x|^{2}\rangle
	&\leq C\langle e^{x},\,|u - w|^{2}\rangle.
\end{align}
For $\varepsilon_1,\varepsilon_2$ suitable, we have
\begin{align}\label{eq60}
	&\frac{d}{dt}\langle e^{x},\,|u - w|^{2}\rangle \leq C \langle e^{x},\,|u - w|^{2}\rangle.
\end{align}
By Gronwall inequality, we obtain
\begin{align*}
	\|u-w\|_{L^2(0,1)}^2\leq \langle e^{x},\,|u - w|^{2}\rangle\leq 0
\end{align*}
then $u=w$.
\renewcommand{\theequation}{\thesection.\arabic{equation}}
\setcounter{equation}{0}
\section{Stability}
\noindent
Let $q\in\mathbb{R}[x]$ given by $q(x)=1+4x-x^3$, now we will prove the exponential decay of the solutions for the $L^{2}$-norm. For this we want to obtain a priori estimate time-independent for this norm.
\begin{theorem}
	If $u$ is solution of \eqref{eq47}-\eqref{eq50}, then
	\begin{align}
		\|u\|_{L^2(0,1)}^2 \leq Ce^{- t}
	\end{align}
	where $C$ depends only of $u_0$.
\end{theorem}
\proof Multiplying \eqref{eq47} by $q(x)\overline{u}$, we obtain
\begin{align}
	&i\int_{0}^{1}q\overline{u}u_{t}dx + i\chi \mathfrak{B}(t)\int_{0}^{1}q\overline{u}u_{xxx}dx + \gamma  \mathfrak{C}(t)\int_{0}^{1}q\overline{u}u_{xx}dx \nonumber \\
	\label{603}& + \int_{0}^{1}I(x,\,t)q\overline{u}u_{x}dx - \int_{0}^{1}q|u|^{4}dx = 0.
\end{align}
Each term is treated separately, using integrating by parts and using \eqref{eq48}-\eqref{eq50}.
\begin{align*}
	i\chi \mathfrak{B}(t)\int_{0}^{1}q\overline{u}u_{xxx}dx&=-i\chi\mathfrak{B}(t)\int_0^1u_{xx}q_x\overline{u} \ dx-i\chi\mathfrak{B}(t)\int_0^1u_{xx}q\overline{u}_x \ dx\\
	&=i\chi\mathfrak{B}(t)\int_0^1u_{x}q_{xx}\overline{u} \ dx+i\chi\mathfrak{B}(t)\int_0^1|u_{x}|^2q_x \ dx
	-i\chi\mathfrak{B}(t)\int_0^1u_{xx}q\overline{u}_x \ dx.\\
	\gamma  \mathfrak{C}(t)\int_{0}^{1}q\overline{u}u_{xx}dx&=-\gamma  \mathfrak{C}(t)\int_{0}^{1}q_x\overline{u}u_{x}dx-\gamma  \mathfrak{C}(t)\int_{0}^{1}q|u_{x}|^2 \ dx.
\end{align*}
Replacing \eqref{603} we have
\begin{align*}
	&i\int_{0}^{1}q\overline{u}u_{t}dx + i\chi \mathfrak{B}(t)\int_{0}^{1}q_{xx}\overline{u}u_{x}dx + i\chi \mathfrak{B}(t)\int_{0}^{1}q_{x}|u_{x}|^{2}dx \nonumber \\
	&- i\chi \mathfrak{B}(t)\int_{0}^{1}q\overline{u}_{x}u_{xx}dx -
	\gamma  \mathfrak{C}(t)\int_{0}^{1}q_{x}\overline{u}u_{x}dx -
	\gamma  \mathfrak{C}(t)\int_{0}^{1}q|u_{x}|^{2}dx \nonumber \\
	& + \int_{0}^{1}I(x,\,t)q\overline{u}u_{x}dx - \int_{0}^{1}q|u|^{4}dx = 0.
\end{align*}
Applying conjugate
\begin{align*}
	&-i\int_{0}^{1}qu\overline{u}_{t}dx - i\chi \mathfrak{B}(t)\int_{0}^{1}q_{xx}u\overline{u}_{x}dx - i\chi \mathfrak{B}(t)\int_{0}^{1}q_{x}|u_{x}|^{2}dx \nonumber \\
	&+ i\chi\mathfrak{B}(t)\int_{0}^{1}qu_{x}\overline{u}_{xx}dx -
	\gamma  \mathfrak{C}(t)\int_{0}^{1}q_{x}u\overline{u}_{x}dx -
	\gamma  \mathfrak{C}(t)\int_{0}^{1}q|u_{x}|^{2}dx \nonumber \\
	& + \int_{0}^{1}\overline{I(x,\,t)}qu\overline{u}_{x}dx - \int_{0}^{1}q|u|^{4}dx = 0.
\end{align*}
Subtracting the above equations yields
\begin{align*}
	&i\frac{d}{dt}\|\sqrt{q}\,u\|_{L^{2}(0,\,1)}^{2} + i\chi \mathfrak{B}(t)\int_{0}^{1}q_{xx}\frac{d}{dx}(|u|^{2})dx + 2i\chi \mathfrak{B}(t)\int_{0}^{1}q_{x}|u_{x}|^{2}dx \nonumber \\
	& - i\chi \mathfrak{B}(t)\int_{0}^{1}q\frac{d}{dx}(|u_{x}|^{2})dx -
	\gamma 2i\mathfrak{C}(t)Im\int_{0}^{1}q_{x}u_{x}\overline{u} \ dx \nonumber \\
	& + 2iIm\int_{0}^{1}I(x,\,t)q\overline{u}u_{x} \ dx = 0.
\end{align*}
Hence, integrating by part some terms we obtain
\begin{align*}
	&\ \frac{d}{dt}\|\sqrt{q}\,u\|_{L^{2}(0,\,1)}^{2} - \chi \mathfrak{B}(t)\int_{0}^{1}q_{xxx}\,|u|^{2}dx + 3\chi \mathfrak{B}(t)\|\sqrt{q_{x}}\,u_{x}\|_{L^{2}(0,\,1)}^{2} \nonumber \\
	& -2\gamma \mathfrak{C}(t)Im\int_{0}^{1}q_{x}\overline{u}u_{x}dx + 2 Im\int_{0}^{1}I(x,\,t)q\overline{u}u_{x}dx =  0
\end{align*}
or
\begin{align*}
	&\frac{d}{dt}\|\sqrt{q}\,u\|_{L^{2}(0,\,1)}^{2} + 6\chi \mathfrak{B}(t)\int_{0}^{1}|u|^{2}dx + 3\chi \mathfrak{B}(t)\|\sqrt{q_{x}}\,u_{x}\|_{L^{2}(0,\,1)}^{2} \nonumber \\
	& = 2\gamma  \mathfrak{C}(t)Im\int_{0}^{1}q_{x}\overline{u}u_{x}dx - 2Im\int_{0}^{1}I(x,\,t)q\overline{u}u_{x}dx.
\end{align*}
Each term in second member is treated separated
\begin{align*}
	&2\gamma  \mathfrak{C}(t)Im\int_{0}^{1}q_{x}\overline{u}u_{x}dx\leq C\int_{0}^{1}|q_{x}\overline{u}u_{x}|dx\leq C \|q_{x}\,u_{x}\|_{L^{2}(0,\,1)} \ \|u\|_{L^{2}(0,\,1)}\\
	&\leq \varepsilon_1\|q_{x}\,u_{x}\|_{L^{2}(0,\,1)}^{2}+C(\varepsilon_1)\|u\|_{L^{2}(0,\,1)}^{2}
	\leq \varepsilon_1\sup_{x\in(0,1)}|q_x|\int_0^1|q_x||u_x|^2 \ dx\\
	& \ \ +C(\varepsilon_1)\|u\|_{L^{2}(0,\,1)}^{2}
	\leq \varepsilon_1C\|\sqrt{q_{x}}\,u_{x}\|_{L^{2}(0,\,1)}^{2}+C(\varepsilon_1)\|u\|_{L^{2}(0,\,1)}^{2}
\end{align*}
and
\begin{align*}
	&- 2Im\int_{0}^{1}I(x,\,t)q\overline{u}u_{x}dx\leq 2\sup_{x\in(0,1)}|q| \ \sup_{Q}|I(x,\,t)|
	\int_{0}^{1}|u| \ |\sqrt{q_{x}}u_{x}|dx\\
	&\leq C\|u\|_{L^{2}(0,\,1)}^{2}\|\sqrt{q_{x}}\,u_{x}\|_{L^{2}(0,\,1)}\leq
	\varepsilon_2\|u\|_{L^{2}(0,\,1)}^{2}+C(\varepsilon_2)\|\sqrt{q_{x}}\,u_{x}\|_{L^{2}(0,\,1)}^{2}.
\end{align*}
Let $k=\inf_{t\in[0,T]}\mathfrak{B}(t)$, then
\begin{align*}
	&\frac{d}{dt}\|\sqrt{q}\,u\|_{L^{2}(0,\,1)}^{2} + 6\chi k\|u\|_{L^2(0,1)}^2 + 3\chi k\|\sqrt{q_{x}}\,u_{x}\|_{L^{2}(0,\,1)}^{2} \nonumber\\
	& \leq \varepsilon_1C\|\sqrt{q_{x}}\,u_{x}\|_{L^{2}(0,\,1)}^{2}+C(\varepsilon_1)\|u\|_{L^{2}(0,\,1)}^{2}
	+\varepsilon_2\|u\|_{L^{2}(0,\,1)}^{2}+C(\varepsilon_2)\|\sqrt{q_{x}}\,u_{x}\|_{L^{2}(0,\,1)}^{2}
\end{align*}
implies que
\begin{align*}
	&\frac{d}{dt}\|\sqrt{q}\,u\|_{L^{2}(0,\,1)}^{2} + (6\chi k-C(\varepsilon_1)-\varepsilon_2)\|u\|_{L^2(0,1)}^2 + (3\chi k -\varepsilon_1C-C(\varepsilon_2))\|\sqrt{q_{x}}\,u_{x}\|_{L^{2}(0,\,1)}^{2} \leq 0.
\end{align*}
Let $C_1=6\chi k-C(\varepsilon_1)-\varepsilon_2$ and $C_2=3\chi k -\varepsilon_1C-C(\varepsilon_2)$. For $\varepsilon_1$ and $\varepsilon_2$ suitable we have $C_1$ and $C_2$ are positive numbers, then, integrating over $(0,t)$ with $t\leq T$
\begin{align*}
	\|\sqrt{q}\,u\|_{L^{2}(0,\,1)}^{2}\leq \|\sqrt{q}\,u_0\|_{L^{2}(0,\,1)}^{2}-C_1\int_0^t\|u\|_{L^{2}(0,\,1)}^{2}.
\end{align*}
Since $\inf_{x\in(0,1)}|q|=1$, we have
$$\|u\|_{L^{2}(0,\,1)}^{2}\leq \|\sqrt{q}\,u\|_{L^{2}(0,\,1)}^{2}$$
then
\begin{align*}
	\|u\|_{L^{2}(0,\,1)}^{2}\leq \|\sqrt{q}\,u_0\|_{L^{2}(0,\,1)}^{2}-C_1\int_0^t\|u\|_{L^{2}(0,\,1)}^{2}
\end{align*}
the result follow.

\section{Numerical Scheme}
\subsection{Notation.}
For the sake of the following analysis, and for a given $M\in\NN$, we will introduce the vector space
\begin{equation*}
	X_M := \big\{u = [u_0\:  u_1 \: \dots \: u_M]^T \in \CC^{M+1} : u_0 = u_{M-1} = u_{M} = 0  \big\}
\end{equation*}

Let us introduce the classical finite differences operators for complex-valued arrays:  
\begin{align*}
	\big[\bfD^+ u\big]_j = D^+u_j &:= \dfrac{u_{j+1} - u_j}{\Delta x} \\
	\big[\bfD^- u\big]_j = D^-u_j &:= \dfrac{u_{j} - u_{j-1}}{\Delta x} \\
	\bfD_x u  &:= \dfrac{1}{2}\Big(\bfD^+u + \bfD^- u \Big) \\
	\bfD_x^2 u &:= \bfD^+\bfD^- u \\
	\bfD_x^3 u  &:= \bfD\bfD^+\bfD^- u.
\end{align*}


\subsection{Foundations of the Numerical Scheme.}

\noindent For a given $\Delta t < 1 $, $t_n = n\Delta t$, and for $u^n \in X_M$, our numerical scheme is inspired in the discretization of
\eqref{301}, and it is defined as follows:
\[i\dfrac{u^{n+1}-u^n}{\Delta t} + i\xi \bfp^3(t_n) \bfD_x^3u^{n+\frac{1}{2}} + \gamma \bfp^2(t_n) \bfD_x^2u^{n+\frac{1}{2}} + i\bfI^{n}\bfD u^{n}-|u^{n+\frac{1}{2}}|^2u^{n+\frac{1}{2}} = i\bfJ^{n+\frac{1}{2}},\]
where $\bfJ^{n} \in X_M: [\bfJ^{n}]_i = J(x_i,t_{n})$, while $\bfI^n \in \RR^{n\times n}$ such that
\[ [\bfI^n]_{ij} = \begin{cases}I(x_i,t_n),\quad \text{ if }i=j \\ 0 \quad \text{ in other case. }  \end{cases} \]
As this is a nonlinear problem, we will find approximate solutions using a Picard Fixed Point iteration. Re-arranging terms and using matrix notation, we get
\begin{equation}\label{esq_num}
	\begin{matrix}
		\ \big[\bfId - 2i\Delta t \xi \bfp^3(t_n) \bfD^3 -2i\Delta t \gamma \bfp^2(t_n) \bfD^2\big]u^{n+1} \\
		= \big[\bfId -2i\Delta t \xi \bfp^3(t_n) \bfD^3 + 2i\Delta t \gamma \bfp^2(t_n)\bfD^2 - 2\Delta t \bfI^n\bfD  \big]u^n  - 2i\Delta t |u^{n+\frac{1}{2}}|^2u^{n+\frac{1}{2}} + 2\Delta t \bfJ^{n+\frac{1}{2}}.
	\end{matrix} 
\end{equation}
where $\bfD$, $\bfD^2$ and $\bfD^3$ are matrices in $\RR^{2M+1\times 2M+1}$ operating over complex-valued vectors, while $\bfI$ is the identity matrix. Due to the boundary conditions, we have
\begin{align*}
	\bfD &= \dfrac{1}{\Delta x}
	\left[
	\begin{matrix}
		0  & 1 & & &  \\
		-1 & 0 & \ddots & \\
		& \ddots & \ddots & \ddots & \\
		&        & \ddots &   0    &  1 \\
		&        &        &   -1   &  0 
	\end{matrix}
	\right] \\
	\bfD^2 &= \dfrac{1}{\Delta x^2}
	\left[
	\begin{matrix}
		-2  & 1 & & &  \\
		1 & -2 & \ddots & \\
		& \ddots & \ddots & \ddots & \\
		&        & \ddots &   -2    &  1 \\
		&        &        &   1   &  -2 
	\end{matrix}
	\right] \\
	\bfD^3 &= \dfrac{1}{\Delta x^3}
	\left[
	\begin{matrix}
		0  & -1 & \frac{1}{2} & & &  &   \\
		1 & 0 & -1 & \frac{1}{2} & & &    \\
		-\frac{1}{2} & 1 & 0 & -1 & \frac{1}{2} & \\
		& \ddots & \ddots & \ddots & \ddots & \ddots & \\
		& & -\frac{1}{2}   & 1 & 0 & -1 & \frac{1}{2} \\
		& & &  -\frac{1}{2} & 1 &   0    &  -1 \\
		&    & &   &  -\frac{1}{2}      &   1   &  0 
	\end{matrix}
	\right] \\
\end{align*}  

To compute the numerical solution, we will use a fixed-point method in order to solve equation \eqref{esq_num} for each time-step. As in Delfour, Fortin and Payre, for a $u^{p=1}= u^n\in X_M$ given, we compute a sequence of complex vectors $\{u^p\},p=2,3,4,\dots,$ until a stopping criteria is verified. The sequence is given by
\begin{equation}\label{pto_fijo}
	\begin{matrix}
		\ \big[\bfId - 2i\Delta t \xi \bfp^3(t_n) \bfD^3_x -2i\Delta t \gamma \bfp^2(t_n) \bfD^2_x\big]u^{p} \\
		= \big[\bfId -2i\Delta t \xi \bfp^3(t_n) \bfD^3_x + 2i\Delta t \gamma \bfp^2(t_n)\bfD^2_x - 2\Delta t \bfI^n\bfD  \big]u^n  - 2i\Delta t \Big|\frac{u^p+u^n}{2}\Big|^2\frac{u^p+u^n}{2} + 2\Delta t \frac{\bfJ^{p}+\bfJ^n}{2}
	\end{matrix}
\end{equation}
In other words, we have to solve a linear system of equations many times per timestep until a stopping criterion is fulfilled, where the matrix to be inverted has a pentadiagonal structure. The fixed point scheme can stop by two reasons: one, if we have
\begin{equation*}
	||u^p - u^{p-1}||_{2}<\hat{\delta}
\end{equation*}
where $\hat{\delta}$ is given. In our computations: $\hat{\delta} = 10^{-14}$. In that case, we do $u^{n+1} = u^p$ and then proceed to the next timestep. The second reason used to stop the iteration is if \eqref{pto_fijo} is not a contraction anymore; that is,
\begin{equation*}
	\frac{||u^p - u^{p-1} ||_2}{ ||u^{p-1} - u^{p-2} ||_2} \geq 1
\end{equation*}
In that case the scheme cannot guarantee the uniqueness of the numerical solution, and thus, the computation ends with no output. The scheme has linear convergence for $\Delta t$ sufficiently small. \\

\noindent The function that transforms $u(x,t)$ to $v(\xi,\tau)$ is used only when we need to plot the numerical solution.

\subsection{Properties}

One of the most remarkable properties of the scheme is that the $L^2-$norm of the solution depends on the length of the domain. We have the following result.
\begin{teo}
	For $u^n \in X_M$ the solution of problem \eqref{esq_num}, there exists a function $C(t)$, depending on the boundaries $\alpha(t)$ and $\beta(t)$, such that
	\[||u^{n+1}||_2^2 \leq C(t)||u^n||_2^2 \]
\end{teo}
{\it Proof:} doing the inner product between \eqref{esq_num} and $u^{n+\frac{1}{2}}$ and extracting the imaginary part, we get
\[||u^{n+1}||_2^2 = ||u^n||_2^2 + Im\Big(i\big(\bfJ^{n+\frac{1}{2}},u^{n+\frac{1}{2}} \big)\Big), \]
where

\begin{align}
	i\big(\bfJ^{n+\frac{1}{2}},u^{n+\frac{1}{2}} \big) &= -(2\gamma \bfp^2(t)\Psi,u^{n+\frac{1}{2}}) - 2i(\bfI\Psi x,u^{n+\frac{1}{2}}) - i(\bfI \Phi,u^{n+\frac{1}{2}}) + i(G'(t),u^{n+\frac{1}{2}}) \\
	& - (2|u^{n+\frac{1}{2}}|^2G(x,t),u^{n+\frac{1}{2}}) - (u^{n+\frac{1}{2}}G,|u^{n+\frac{1}{2}}|^2) \\
	&-2(G^2,|u^{n+\frac{1}{2}}|^2) - (G^2,(u^{n+\frac{1}{2}})^2) - (G^3,u^{n+\frac{1}{2}}).
\end{align}

\noindent Our objective is to obtain an upper bound for $G$, $I$, $\Psi$ and $\Phi$ terms. We have
\[\phi(t) = \varphi_3(t) [\bfp(t)]^{-1} = u_x(1,t)\bfp(t)[\bfp(t)]^{-1} = u_x(1,t) \]
and thus, we can bound $\phi(t)$ assuming that $u_x(1,t)$ is also bounded. Furthermore,

\[ |\Phi(t)| = |2(\varphi_1 - \varphi_2) + \phi| = |2(u(0,t)-u(1,t)) + u_x(1,t)| \leq C_\Phi \]
\[ |\Psi(t)| = |\varphi_2 - \varphi_1 - \phi| = |u(1,t)-u(0,t) - u_x(1,t)| \leq C_\Psi \]
this is, the $G(x,t)$ function will be bounded if the boundary conditions also are. Thus, we can infer the existence of a real function $G(x,t)$, depending on the boundary conditions and the boundary itself, such that
\[|G(x,t)| \leq C_G(x,t). \]
$G'(t)$ will also be bounded if the same conditions fulfill along with the fact that boundary conditions are in $C^2([9,+\infty])$. Meanwhile,
\[I(x,t) = \dfrac{\bfp(t)}{\bfp'(t)}x - \alpha'(t)\bfp(t) \leq \bfp(t)\Big|\dfrac{1}{\bfp'(t)} - \alpha'(t) \Big|. \]
Through Cauchy-Schwarz, we have the following estimates

\begin{align}
	\big|( \bfp^2(t)\Psi,u^{n+\frac{1}{2}})\big| &\leq  \bfp^2(t)||u^{n+\frac{1}{2}}||_2 C_\Psi \\
	|(\bfI\Psi x,u^{n+\frac{1}{2}})|  &\leq \bfp^2(t)\Big|\dfrac{1}{\bfp'(t) - \alpha'(t)} \Big|||u^{n+\frac{1}{2}}||_2 C_\Psi \\
	|(\bfI \varphi,u^{n+\frac{1}{2}})|  &\leq \bfp(t)\Big|\dfrac{1}{\bfp'(t) - \alpha'(t)} \Big|||u^{n+\frac{1}{2}}||_2 |\varphi_1(t)| \\
	|(G'(t),u^{n+\frac{1}{2}}| &\leq |G'(t)| ||u^{n+\frac{1}{2}}||_1  \\
	|(|u^{n+\frac{1}{2}}|^2G(x,t),u^{n+\frac{1}{2}})| &\leq C_G(x,t) ||u^{n+\frac{1}{2}}||_3^3   \\
	(u^{n+\frac{1}{2}}G(x,t),|u^{n+\frac{1}{2}}|^2) &\leq C_G(x,t) ||u^{n+\frac{1}{2}}||_3^3  \\
	Im(G^2,|u^{n+\frac{1}{2}}|^2) &= 0 \\
	(G^2,(u^{n+\frac{1}{2}})^2) & \leq C_G(x,t)^2 ||u^{n+\frac{1}{2}}||^2_2  \\
	(G^3,u^{n+\frac{1}{2}}) & \leq C_G(x,t)^3 ||u^{n+\frac{1}{2}}||_1 
\end{align}

Because $||z||_3^3 \leq ||z||_2^2 ||z||_\infty,\: z \in X_M$, we can thus infer the existence of a function $C(t)$, bounding both $C_G(x,t)$ and the boundary conditions, such that we get the desired result. 

\subsection{Convergence.}
Let us recall that $u_j^n \approx u(x_j,t_n)$, for $u(x,t)$ the exact solution of our problem. We will denote $\bfu(x,t)\in \RR^M:[\bfu(x,t)]_k := u(x_k,t)$. We will define the truncation error $\bfF^n\in \RR^M$ as follows
\begin{equation}
	\begin{matrix}
		\ \big[\bfI - 2i\Delta t \xi \bfp^3(t_n) \bfD^3 -2i\Delta t \gamma \bfp^2(t_n) \bfD^2\big]\bfu (x,t_{n+1}) \\
		= \big[\bfI -2i\Delta t \xi \bfp^3(t_n) \bfD^3 + 2i\Delta t \gamma \bfp^2(t_n)\bfD^2 - 2\Delta t \bfI^n\bfD  \big]\bfu (x,t_n)\\  - 2i\Delta t |\bfu(x,t_{n+\frac{1}{2}})|^2\bfu(x,t_{n+\frac{1}{2}}) + 2\Delta t \bfJ^{n+\frac{1}{2}} + \bfF^n.
	\end{matrix} 
\end{equation}
The following result holds.
\begin{teo}
	For $u_j^n$ the solution of problem \eqref{esq_num}, and $u(x,t)$ the solution of problem (...), such that $u \in C^1[(0,T),C^2(\Omega)]$. Then, there exists a constant $C$, depending on $\Omega$, such that
	\[ || \bfF^n||_2 \leq C\O(\Delta t + \Delta x^2). \]
\end{teo}
{\it Proof:} due to Taylor expansions at $(x_k,t_n)$, it is well known that
\[D_tu(x_k,t_n) = u_t(x_k,t_n) + \O(\Delta t),\quad D^2_xu(x_k,t_n)= u_{xx}(x_k,t_n) + \O(\Delta x^2) \]
\[D^3_xu(x_k,t_n) = u_{xxx}(x_k,t_n) + \O(\Delta x^2) ,\quad u(x_k,t_{n+\frac{1}{2}}) = u(x_k,t_n) + \O(\Delta t).\]
Thus, at $x_k = x_0 + k\Delta x$ we have
\begin{equation*}
	\begin{matrix}
		iu_t(x_k,t_n) +i\xi \bfp^3(t_n) u_{xxx}(x_k,t_n) + \gamma \bfp^2(t_n)u_{xx}(x_k,t_n) + i \bfI^nu_x(x_k,t_n) - |u(x_k,t_{n+\frac{1}{2}})|^2u(x_k,t_{n+\frac{1}{2}})  \\
		=  i\bfJ^{n+\frac{1}{2}} + F_k^n + \O(\Delta t + \Delta x^2) 
	\end{matrix} 
\end{equation*}
Thus, we get
\[|F_k^n| = \O(\Delta t + \Delta x^2) \]
multiplying by $\Delta x$ and summing for $k$ we get the result.

\subsection{Numerical experiments.}
\begin{figure}
	\centering
	\includegraphics[scale=0.6]{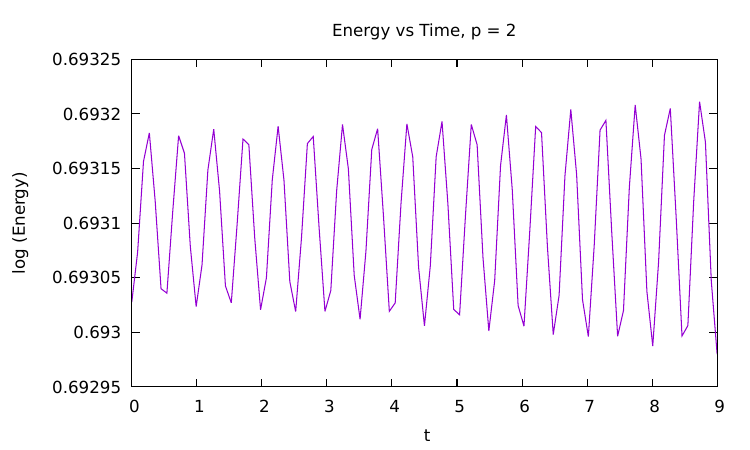}
	\includegraphics[scale=0.6]{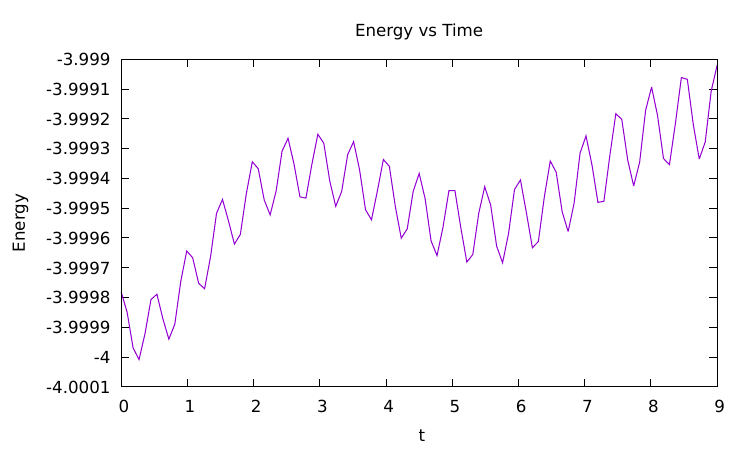}
	\caption{Left: Energy at $L^2$ level. Right: Energy at $H$ level.}
	\label{fig1}
\end{figure}

\begin{figure}
	\centering
	\includegraphics[scale=1]{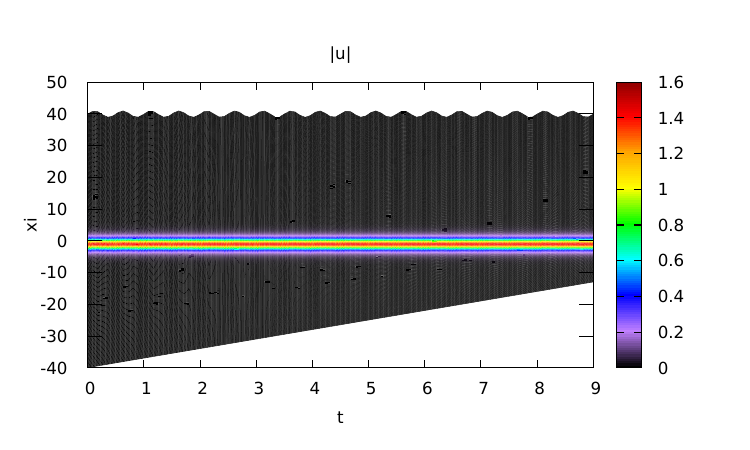}
	\caption{Numerical solution.}
	\label{fig2}
\end{figure}

\noindent {\bf Case 1: travelling soliton.} We have applied our scheme to the problem
\[
\begin{aligned}
iv_\tau -v_{\xi\xi} &= |v^2|v,\qquad (\xi,\tau) \in [-40,40]\times (0,9] \\
v(\xi,0) &= \sqrt{2}\sech(\xi +1) \\
v(\alpha(\tau)),\tau) &= v(\beta(\tau),\tau) = v_\xi(\beta(\tau),\tau) = 0
\end{aligned}
\]
for $\alpha(\tau) = 3\tau - 40$ and $\beta(\tau) = 40 + \sin(4\pi\tau)$. Figure \ref{fig1} illustrates the situation of the energy at $L^2$ and $H^1$ levels, respectively. Figure \ref{fig2} shows the numerical solution obtained.

\noindent {\bf Case 2: travelling soliton for a returning linear boundary condition.} We will now consider the case
\[
\begin{aligned}
iv_\tau +v_{\xi\xi} &= |v^2|v,\quad (\xi,\tau)\in (\alpha(\tau),\beta(\tau)) \times (0,T] \\
v(\xi,0) &= 2\sqrt{2}e^{i\xi-5}\sech(2\xi +5) \\
v(\alpha(\tau)),\tau) &= v(\beta(\tau),\tau) = v_\xi(\beta(\tau),\tau) = 0
\end{aligned}
\]
for $\alpha(\tau) = \dfrac{x_0+20}{1.7}|\tau-1.7|-20$ and $\beta(\tau) = x_f + \sin(4\pi\tau)$. This is an interesting case, as the $L^2-$norm is not preserved during the whole integration time, as one would expect. Figure \ref{fig2} shows our results. We can clearly see that the $L^2-$ norm is lost when the soliton approaches to the narrow zone, but then it goes back when the domain recovers its original length.

\begin{figure}[htbp]
	\centering
	\includegraphics[scale=0.6]{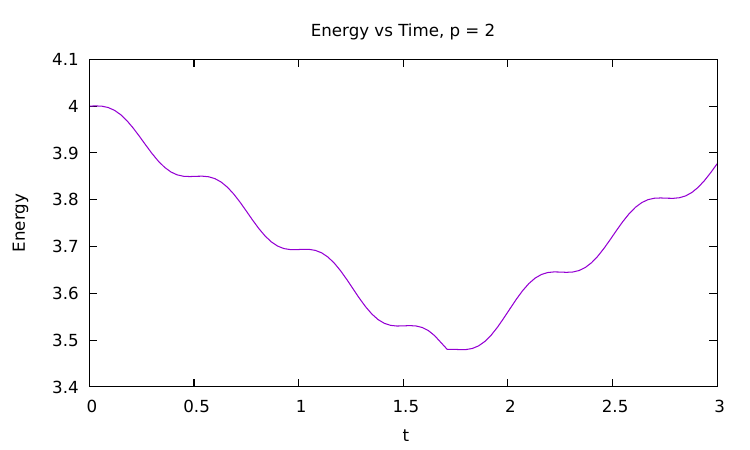} \\
	\includegraphics[scale=1.5]{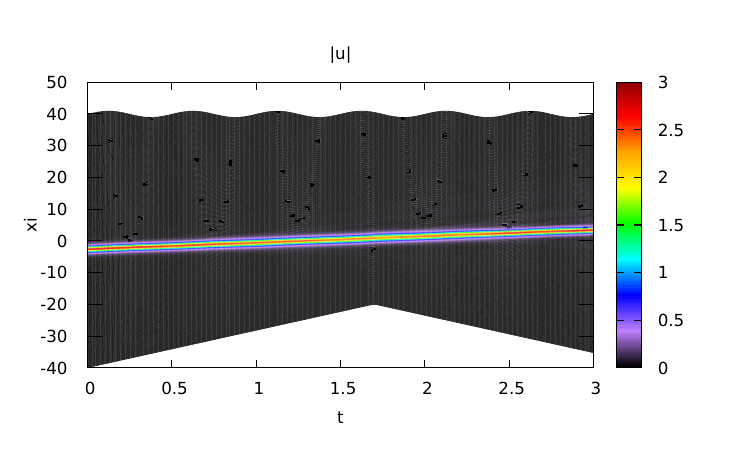}
	\caption{Left: Energy at $L^2$ level. Right: Numerical solution. In both cases, $x_0 = -40$ and $x_f = 40$.}
	\label{fig3}
\end{figure}

Table \ref{t1} illustrates how the difference in $L^2-$norm evolves with respect to the initial length of the domain. We can clearly see that the difference is small when the initial length is larger, which indicates a dependence in the preservation of the $L^2-$ norm on the length of the domain.

\begin{center}
	\begin{table}
		\begin{tabular}{|c|c|c|c|}
			\hline
			$x_0$ & $x_f$ & $x_f-x_0$ & $\Delta || u ||_2$ \\ \hline
			$-40$ & $10$ & $50$  & $0.87353$ \\ \hline
			$-40$ & $20$ & $60$  & $0.71152$ \\ \hline
			$-40$ & $40$ & $80$  & $0.52033$ \\ \hline
			$-40$ & $60$ & $100$  & $0.41055$  \\ \hline
		\end{tabular}
		\caption{Difference in $L^2-$ norm of the numerical solution for differente values of $x_f$ for the second case.}
		\label{t1}
	\end{table}
\end{center}

\noindent {\bf Case 3: two soliton solution.} As a final case, we will repeat the previous case but for the initial condition
\[  v(\xi,0) = 2\sqrt{2}e^{i\xi-5}\sech(2\xi +5) + 2\sqrt{2}e^{i\xi+10}\sech(2\xi-10) \]

\noindent Figure \ref{fig3} shows the numerical solution obtained by the scheme. While we can see more detail as the domain gets narrower, the $L^2-$ norm gets smaller. 

\begin{figure}[htbp]
	\centering
	\includegraphics[scale=0.6]{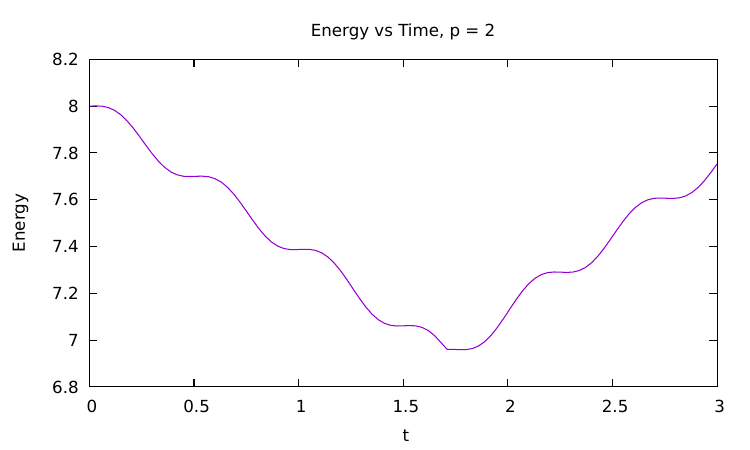} \\
	\includegraphics[scale=1.5]{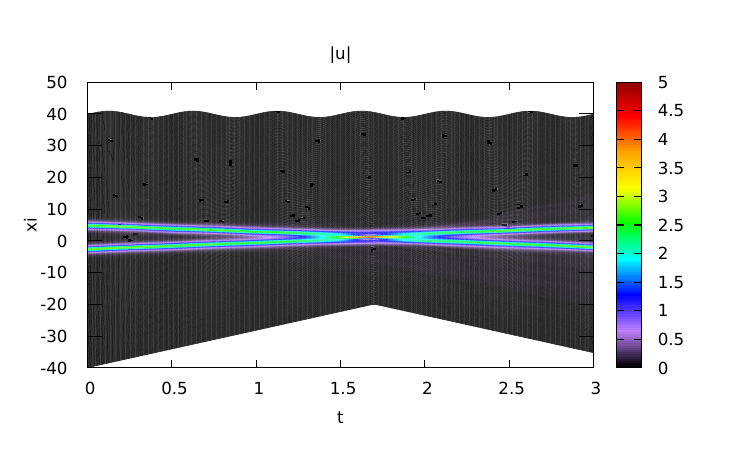}
	\caption{Numerical solution for the two-soliton case, with its respective $L^2-$ time evolution plot.}
	\label{fig4}
\end{figure}

\newpage

\begin{appendices}
\section{Appendix: The Gauge transformation}

We consider the Gauge transformation
\begin{align}
\label{108}v(\xi,\,\tau) = e^{id_{2}\xi + id_{3}\tau}w(\xi - d_{1}\tau,\,\tau)\equiv e^{\Theta}w(y,\,s)
\end{align}
where $\Theta = id_{2}\xi + id_{3}\tau,$\ $y = \xi - d_{1}\tau$ and $s = \tau.$ Then
\begin{align*}
& v_{\tau} = id_{3}e^{\Theta}w - d_{1}e^{\Theta}w_{y} + e^{\Theta}w_{s};
\ v_{\xi} = id_{2}e^{\Theta}w + e^{\Theta}w_{y}; \\
&v_{\xi\xi} = -d_{2}^{2}e^{\Theta}w + 2id_{2}e^{\Theta}w_{y} +e^{\Theta}w_{yy};
\\ 
& v_{\xi\xi\xi} = -id_{2}^{3}e^{\Theta}w - 
3d_{2}^2e^{\Theta}w_{y}
+3id_{2}e^{\Theta}w_{yy} + e^{\Theta}w_{yyy}.
\end{align*}
Replacing in \eqref{101} we obtain
\begin{align*}
-d_{3}e^{\Theta}w - id_{1}e^{\Theta}w_{y} + ie^{\Theta}w_{s} - \gamma d_{2}^{2}e^{\Theta}w + 2i\gamma d_{2}e^{\Theta}w_{y} + \gamma e^{\Theta}w_{yy}
\MoveEqLeft[6]
\\
+ \chi d_{2}^{3}e^{\Theta}w_{y} 
- 3i\chi d_{2}^{2}e^{\Theta}w_{y} 
- 3\chi d_{2}e^{\Theta}w_{yy} 
+ i\chi e^{\Theta}w_{yyy} 
\ = \ |z|^{2}e^{\Theta}z.
\end{align*}
Gathering similar terms and performing straightforward calculation we obtain
\begin{align*}
d_{1} = \frac{\gamma^{2}}{3\chi}\ ;\ d_{2} = \frac{\gamma}{3\chi}\ ;\ d_{3} = \frac{-2\gamma^{3}}{27\chi^{2}}.
\end{align*}
Thereby in \eqref{101} follows that
\begin{align*}
w_{s} + \chi w_{yyy}  + i|w|^{2}w = 0,
\end{align*}
This way using the Gauge transformation we have the equivalent problem to \eqref{101}
\begin{align*}
\begin{cases}
w_{s} + \chi w_{yyy}  + i|w|^{2}w = 0,
\quad \mbox{in } \widetilde{Q}_{s},   \\
w(y,0) = e^{-i\frac{\gamma}{3}y}v_{0}(y)\quad\mbox{in }
\widetilde{I}_{0}, \\
v(\widetilde{\alpha}(s),\,s) = \varphi_{1}(s),\ v(\widetilde{\beta}(s),\,s) = \varphi_{2}(s),\ v_{y}(\widetilde{\beta}(s),\,s) = \varphi_{3}(s).
\end{cases}
\end{align*}
where $\widetilde{\alpha}(s)=\alpha(s)-d_1s$,
$\widetilde{\beta}(s)=\beta(s)-d_1s$,
$
\widetilde{I}_{s} = \{y\in\mathbb{R}:\,
\widetilde{\alpha}(s) < y < \widetilde{\beta}(s)\}$ and 
$\widetilde{Q}_{s} = \{(y,\,s)\in\mathbb{R}^{2}:y\in
\widetilde{I}_{s},\,0 < s < T,\,T>0\}.$ 
The above Gauge transformation is a bi-continuous map from $L^{p}([0,\,T]: H^{\sigma}(Q_\tau)$ to $L^{p}([0,\,T]: H^{\sigma}(\widetilde{Q}_{s})$ , as long as $0<T<+\infty.$ 
Using the Gauge transformation we see that the 
$v_{\xi\xi}$ term in the equation \eqref{101} can be removed,
Thus, using the transform moving boundary of the chapter 2,
the system \eqref{eq47}-\eqref{eq50} is reduced to
\begin{align*}
\begin{cases}
u_{t} + \chi[{\widetilde{\bf p}}(t)]^{3}u_{xxx} -i \widetilde{I}(x,\,t)u_{x} = |u|^{2}u + i\widetilde{J}(x,\,t),
& 0<x<1, \quad t>0\\
u(x,0)=u_0(x),& 0<x<1,\\
u(0,\,t) = u(1,\,t) = u_{x}(1,\,t) = 0
& 
t>0.
\end{cases}
\end{align*}
where $\alpha(\tau)$ and $\beta(\tau)$ are replaced by 
$\widetilde{\alpha}(\tau)$ and $\widetilde{\beta}(\tau)$,
respectively, and $\widetilde{\bf p}(t)$, $\widetilde{I}(x,\,t)$,
$\widetilde{J}(x,\,t)$, are defined in term of these modified 
moving boundary functions.
\end{appendices}


\begin{thebibliography}{10}

\bibitem{alves} 
Alves, M., Sep\'ulveda, M., Vera, O., 
{\it Smoothing properties for the higher-order nonlinear Schr\"odinger equation with constant coefficients}. Nonlinear Anal. 71 (2009), no. 3-4, 948-966.

\bibitem{au} Aubin, J.P. {\it Un theoréme de compacité, C. R. Acad. Sci.} Paris. 256(1963)5042-5044.

\bibitem{ben} Ben-Artzi, M.H., Koch, H., and Saut, J.-C. {\it Dispersion estimates for fourth order Schr\"{o}dinger equations.} C. R. Acad. Sci. Paris S\'er. I Math. 330(2000)87-92.

\bibitem{bi01} Bisognin, E., Bisognin, V., and Vera, O. {\it Stabilization of solutions for the higher order nonlinear Schr\"{o}dinger equation with localized damping.} EJDE. 6(2007)1-18.

\bibitem{bi02} V. Bisognin, V., and Vera, O. {\it On the unique continuation property for the higher order nonlinear Schr\"{o}dinger equation with constant coefficients.}  Turkish J. Math. 30(2007)265-302.

\bibitem{bi} Bisognin, E., Bisognin, V., Sep\'{u}lveda, M., and Vera, O. {\it Coupled system of Korteweg de Vries equations type in domains with moving boundaries.} J. Comput. Appl. Math. 220 (2008), 290-321.

\bibitem{bi1} Bisognin, V., Buriol, C., Ferreira, M.V., Sep\'{u}lveda, M., and Vera, O. {\it Asymptotic behaviour for a nonlinear Schrodinger equation in domains with moving boundaries.} Acta Appl. Math. 125 (2013), 159--172.

\bibitem{bona} Bona, J., and J.C. Saut. {\it Dispersive blow-up of solutions of generalized Korteweg–de Vries equation.} J. Differential Equations. 103(1993)3-57.

\bibitem{xavier} Carvajal, X. and Linares, F. {\it A higher order nonlinear Schr\"{o}dinger equation with variable coefficients.} Differential Integral Equations. 16(2003)1111-1130.

\bibitem{rodrigo1}
Cavalcanti, M.-M., Corr\^ea, W.-J., Faminskii, A.-V., Sep\'ulveda, M.-A., V\'ejar-Asem, R, 
{\it Well-posedness and asymptotic behavior of a generalized higher order nonlinear Schr\"odinger equation with localized dissipation}. 
Comput. Math. Appl. 96 (2021), 188–208.

\bibitem{rodrigo2}
Cavalcanti, M.-M., Corr\^ea, W.-J., 
Sep\'ulveda, M.-A., V\'ejar-Asem, R., 
{\it Finite difference scheme for a higher order nonlinear Schrödinger equation.} 
Calcolo 56 (2019), no. 4, Paper No. 40, 32 pp.

\bibitem{rodrigo3}
Cavalcanti, M.-M., Corr\^ea, W.-J., 
Sep\'ulveda, M.-A., V\'ejar-Asem, R., 
{\it Finite difference scheme for a high order nonlinear Schrödinger equation with localized damping}. 
Stud. Univ. Babe\c s-Bolyai Math. 64 (2019), no. 2, 161–172.

\bibitem{Chen}
Chen, S.-S., Tian, B. Qu, Q.-X. Li, H. 
Sun, Y., Du, X.-X.,
Alfvén solitons and generalized Darboux transformation for a variable-coefficient derivative nonlinear Schr\"odinger equation in an inhomogeneous plasma,
Chaos, Solitons \& Fractals,
Volume 148,
(2021),
111029.

\bibitem{co} Coddington, E.A., and Levinson, N. {\it Theory of Ordinary Differential
Equations.} McGraw-Hill, New York (1955).

\bibitem{colin} Colin, T., and Gisclon, M. {\it An initial-boundary-value problem that approximate the quarter-plane problem for the Korteweg–de Vries equation.} Nonlinear Anal. 46(2001)869-892.

\bibitem{do} Doronin, G. and Larkin, N. {\it KdV equation in domains with moving boundaries.} J. Math. Anal. Appl. 328(2007)503-517.

\bibitem{dys} Dysthe, K.B. {\it Note on a modification to the nonlinear Schrodinger equation for application to deep water waves.} Proc. R. Soc. Lond. Ser. A 369 (1979), 105--114.

\bibitem{fri} Friedman, A. {\it Partial Differential Equations.}
Holt-Rinehart and Winston, New York (1969).

\bibitem{gear1}Gear, J.A. {\it Strong interactions between solitary waves belonging to different wave modes.} Stud. Appl. Math. 72(1985)95-124.

\bibitem{gear2} Gear, J.A., and Grimshaw, R. {\it Weak and strong interactions between internal solitary waves.} Stud. Appl. Math. 70(1984)235-258.

\bibitem{hasegawa}  Hasegawa, A. and Kodama, Y. {\it Nonlinear pulse propagation in a monomode dielectric guide.} IEEE. J. Quant. Elect. 23(1987)510-524.

\bibitem{hiro} Hirota, R. {\it Direct Methods in Soliton Theory.} Springer, Berlin (1980).

\bibitem{ka1} Karpman, V. {\it Stabilization of soliton instabilities by higher order dispersion: fourth-order nonlinear Schr\"{o}dinger-type equations.}  Phys. Rev. E. Vol. 53. 2(1996)1336-1339.

\bibitem{ka2}  Karpman, V., and Shagalov, A. {\it Stability of soliton described by nonlinear Schr\"{o}dinger-type equations with higher-order dispersion.} Physica D. 144(2000)194-210.

\bibitem{kato} Kato, T. {\it On the Cauchy problem for the (generalized) Korteweg–de Vries equations.} Adv. in Math. Suppl. Stud., Stud. Appl. Math. 8 (1983)93-128.

\bibitem{kodama} Kodama, Y. {\it Optical solitons in a monomode fiber.} J. Phys. Stat. 39(1985)597-614.

\bibitem{la} Larkin, N.A. {\it Korteweg–de Vries and Kuramoto–Sivashinsky equations in bounded domains.} J. Math. Anal. Appl. 297(2004)169-185.

\bibitem{laurey} C. Laurey. {\it Le probleme de Cauchy pour une \'{e}quation de Schr\"{o}dinger non-lin\'{e}aire de ordre 3.} C. R. Acad. Sci. Paris 315(1992)165-168.

\bibitem{lions} Lions, J.L. {\it Quelques M\'{e}thodes de R\'{e}solution
des Problemes aux Limites Non Lin\'{e}aires.} Gauthiers-Villars,
Paris (1969).

\bibitem{gronwall} Pachpatte, B.G. {\it Inequalities for differential and Integral
Equations.} Mathematics in Sciences and Engineering, Vol. 197,
Academic Press (1997).

\bibitem{mauri} Sep\'{u}lveda, M. and Vera, O. {\it Numerical methods for a transport equation perturbed by dispersive terms of 3rd and 5th order, Sci. Ser. A: Mathematics} 
(New Series) 13(2006)13-21.

\bibitem{Staffilani} G. Staffilani, G. {\it On the generalized Korteweg–de Vries type equation.} Differential Integral Equations. Vol. 10. 4(1997)777-796.

\bibitem{teman} Teman, R. {\it Sur um probl\'{e}me non
lin\'{e}aire.} J. Math. Pures Appl. 48(1969)159-172.

\bibitem{octavio}
Vera Villagran, O.-P.,
{\it Stabilization for a fourth order nonlinear Schrödinger equation in domains with moving boundaries}. 
J. Partial Differ. Equ. 34 (2021), no. 3, 268-283.

\bibitem{wloka}  Wloka, J. {\it Partielle Differentialgleichungen.} BSB, B.G. Teubner Veralagsgesellschaft. Leipzig. 1982.

\end{thebibliography}
\end{document}